\documentclass[12pt]{article}
\usepackage{amsmath}
\usepackage{amssymb}
\usepackage{vatola}

\def\q{\quad}
\def\qq{\qquad}
\def\qtq#1{\q\t{#1}\q}
\def\mod#1{\ (\text{\rm mod}\ #1)}
\def\t{\text}
\def\f{\frac}
\def\e{\equiv}
\def\b{\binom}
\def\sls#1#2{(\f{#1}{#2})}
 
\def\Ls#1#2{\Big(\f{#1}{#2}\Big)}
\def\ap{\langle a\rangle_p}
\def\bp{\langle b\rangle_p}
\def\qp#1{q_p(#1)}
\let \pro=\proclaim
\let \endpro=\endproclaim

\begin{document}
 \centerline {\bf
Super congruences concerning binomial coefficients and Ap\'ery-like
numbers}
\par\q\newline
\centerline{Zhi-Hong Sun}\newline \centerline{School of Mathematics
and Statistics}\centerline{Huaiyin Normal University}
\centerline{Huaian, Jiangsu 223300, P.R. China} \centerline{Email:
zhsun@hytc.edu.cn} \centerline{Homepage:
http://maths.hytc.edu.cn/szh1.htm}
 \abstract{Let $p$ be a prime with $p>3$, and let $a,b$ be two
 rational $p-$integers. In this paper we
 present general congruences for $\sum_{k=0}^{p-1}\binom
ak\binom{-1-a}k\frac p{k+b}\mod {p^2}$.  For $n=0,1,2,\ldots$ let
$D_n$ and $b_n$ be Domb and Almkvist-Zudilin
 numbers, respectively. We also establish congruences
 for
 $$\sum_{n=0}^{p-1}\frac{D_n}{16^n},\quad \sum_{n=0}^{p-1}\frac{D_n}{4^n},
 \quad
\sum_{n=0}^{p-1}\frac{b_n}{(-3)^n},\quad
\sum_{n=0}^{p-1}\frac{b_n}{(-27)^n}\pmod {p^2}$$ in terms of certain
binary quadratic forms.
 \par\q
 }
 \endabstract
\section*{1. Introduction}
\par Let $p>3$ be a prime. In 2003, based on his work concerning
hypergeometric functions and Calabi-Yau manifolds,
Rodriguez-Villegas [RV] conjectured the following congruences:
$$\aligned &\sum_{k=0}^{p-1}\f{\b{2k}k^2}{16^k}\e
\Ls{-1}p\mod{p^2},\q\sum_{k=0}^{p-1} \f{\b{2k}k\b{3k}k}{27^k}\e\Ls
{-3}p\mod{p^2}\\&\sum_{k=0}^{p-1}\f{\b{2k}k\b{4k}{2k}}{64^k}\e
\Ls{-2}p\mod{p^2},\q \sum_{k=0}^{p-1} \f{\b{3k}k\b{6k}{3k}}{432^k}\e
\Ls{-1}p\mod{p^2},\endaligned\tag 1.1$$ where $\sls ap$ is the
Legendre symbol. These congruences were later confirmed by Mortenson
[M1-M2] via the Gross-Koblitz formula. It is easily seen that (see
[S7])
$$\aligned &\b{-\f 12}k^2=\f{\b{2k}k^2}{16^k},\q\b{-\f 13}k\b{-\f 23}k=
\f{\b{2k}k\b{3k}k}{27^k},
\\&\b{-\f 14}k\b{-\f 34}k=\f{\b{2k}k\b{4k}{2k}}{64^k},
\q \b{-\f 16}k\b{-\f
56}k=\f{\b{3k}k\b{6k}{3k}}{432^k}.\endaligned\tag 1.2$$
 Let $\Bbb
Z$ be the set of integers.  For a prime $p$ let $\Bbb Z_p$ be the
set of rational numbers whose denominator is not divisible by $p$.
In [S7] the author generalized (1.1) by proving that for $a\in\Bbb
Z_p$,
$$\sum_{k=0}^{p-1}\b ak\b{-1-a}k\e (-1)^{\ap}\mod {p^2},\tag 1.3$$
where $\ap\in\{0,1,\ldots,p-1\}$ is given by $a\e \ap\mod p$. In
[S9] the author established the congruence
$$\sum_{k=0}^{p-1}\b ak\b {-1-a}k\f{2a+1}{2k+1}\e 1+2\f{a-\ap}p\mod
{p^2}.\tag 1.4$$

\par Let $p$ be an odd prime and $a,b\in\Bbb Z_p$ with $ab\not\e 0\mod p$.
 In Section 2 we show that for $\bp\le \ap$ and $\bp\le p-1-\ap$,
$$\sum_{k=0}^{p-1}\b
ak\b{-1-a}k\f 1{k+b}\e
\f{p(s+t+1)(s-t)}{b^2(s+1)\b{\ap}{\bp}\b{p-1-\ap}{\bp}} \mod
{p^2},\tag 1.5$$ where $t=(a-\ap)/p$ and $s=(b-\bp)/p\not\e -1\mod
p$. In Section 3 we present a congruence for $\sum_{k=0}^{p-1}\b
ak\b{-1-a}k\f p{k-a}\mod {p^2}$ under the condition $1\le \ap<\f
p2$. In Section 4, using a combinatorial identity we obtain a
congruence for $\sum_{k=0}^{p-1}\b{2k}k^2\f p{16^k(k+b)}\mod {p^2}$
under the condition $\bp>\f p2$. In particular, we prove that
$$\align &\sum_{k=0}^{p-1}\b{2k}k^2\f p{16^k(6k+1)}
\e 4x^2-2p\mod{p^2}\q\t{for}\q p=x^2+3y^2\e 1\mod 3,\\&
\big(2-(-1)^{\f{p-1}2}\big)\sum_{k=0}^{p-1}\b{2k}k^2\f p{16^k(8k+1)}
\e 4x^2-2p\mod{p^2}\ \t{for}\  p=x^2+2y^2\e 1,3\mod 8,
\\&\Big(3\Ls p3-2\Big)\sum_{k=0}^{p-1}\b{2k}k^2\f
p{16^k(12k+1)}\e 4x^2-2p\mod{p^2}\ \t{for}\  p=x^2+4y^2\e 1\mod 4,
\endalign$$
where $x$ and $y$ are integers. We mention that in [S11, Theorem 2.8
and (6)] the author proved that
$$\sum_{k=0}^{p-1}\b{2k}k^2\f p{16^k(4k+1)}\e
4x^2-2p-\f{p^2}{4x^2}\mod {p^3}\q\t{for prime $p=x^2+4y^2\e 1\mod
4$.}$$

\par  In Section 5 we use congruences in Sections 2-4 to prove super
congruences for three Ap\'ery-like numbers.  Let $\{D_n\},\{b_n\}$
and $\{W_n\}$ be Ap\'ery-like numbers given by
$$\align
&D_n=\sum_{k=0}^n\b nk^2\b{2k}k\b{2n-2k}{n-k},\tag 1.6 \\&
b_n=\sum_{k=0}^{[n/3]}\b{2k}k\b{3k}k\b n{3k}\b{n+k}k(-3)^{n-3k},\tag
1.7
\\&W_n=\sum_{k=0}^{[n/3]}\b{2k}k\b{3k}k\b n{3k}(-3)^{n-3k},\tag 1.8
\endalign$$
where  $[x]$ is the greatest integer not exceeding $x$.
 The numbers $D_n$ $(n=0,1,2,\ldots)$ are called Domb numbers, and $b_n$
 $(n=0,1,2,\ldots)$ are called
Almkvist-Zudilin numbers. For $D_n,b_n,W_n$ See A002895, A125143 and
A291898 in Sloane's database ``The On-Line Encyclopedia of Integer
Sequences", and related papers [AT], [CZ], [S8] and [S12].
\par   For positive integers $a,b$ and $n$, if $n=ax^2+by^2$ for some integers
$x$ and $y$, we briefly write that $n=ax^2+by^2$. Let $p$ be a prime
with $p>3$. In [Su2], the author's brother Z.W. Sun posed the
following conjecture:
$$\sum_{n=0}^{p-1}\f{D_n}{16^n}\e \sum_{n=0}^{p-1}\f{D_n}{4^n}
\e \cases 4x^2-2p\mod{p^2}&\t{if $3\mid p-1$ and so
$p=x^2+3y^2$,}\\0\mod {p^2}&\t{if $3\mid p-2$.}
\endcases\tag 1.9$$
 In [S12, Conjecture 4.15]  the author conjectured
that
$$\sum_{n=0}^{p-1}\f{D_n}{4^n}\e\sum_{n=0}^{p-1}\f{D_n}{16^n}\e
4x^2-2p-\f{p^2}{4x^2}\mod {p^3}\qtq{for}p=x^2+3y^2\e 1\mod 3,$$ and
$$\sum_{n=0}^{p-1}\f{D_n}{4^n}\e -2\sum_{n=0}^{p-1}\f{D_n}{16^n}
\e\f{p^2}2\b{\f{p-1}2}{\f{p-5}6}^{-2}\mod {p^3}\qtq{for}p\e 2\mod
3.$$  In [S12, Conjecture 4.16]  the author also conjectured that
$$\align \sum_{n=0}^{p-1}\f{b_n}{(-3)^n}\e
\sum_{n=0}^{p-1}\f{b_n}{(-27)^n}\e \cases 4x^2-2p-\f{p^2}{4x^2}\mod
{p^3}&\t{if $12\mid p-1$ and $p=x^2+9y^2$,}
\\2p-2x^2+\f{p^2}{2x^2}\mod {p^3}&\t{if $12\mid p-5$ and
$2p=x^2+9y^2$}\endcases\endalign$$ and
$$ \sum_{n=0}^{p-1}\f{b_n}{(-3)^n}\e -15
\sum_{n=0}^{p-1}\f{b_n}{(-27)^n} \e\cases -\f
53p^2\b{[p/3]}{[p/12]}^{-2}\mod {p^3}&\t{if $12\mid p-7$,}
\\\f 56p^2\b{[p/3]}{[p/12]}^{-2}\mod {p^3}&\t{if $12\mid p-11$.}\endcases$$
In Section 5 we prove (1.9) and
$$\aligned\sum_{n=0}^{p-1}\f{b_n}{(-3)^n}&\e
\sum_{n=0}^{p-1}\f{b_n}{(-27)^n}\e \sum_{k=0}^{p-1}
\f{\b{2k}k\b{3k}k}{27^k}\cdot\f p{4k+1}
\\&\e \cases 4x^2-2p\mod {p^2}&\t{if $12\mid p-1$ and so
$p=x^2+9y^2$,}
\\2p-2x^2\mod {p^2}&\t{if $12\mid p-5$ and so
$2p=x^2+9y^2$,}\\0\mod {p^2}&\t{if $p\e 3\mod 4$.}
\endcases\endaligned\tag 1.10$$
We also determine $\sum_{n=0}^{p-1}\f{W_n}{(-3)^n}$ modulo $p^2$. In
particular, we show that for $p\e 1\mod 3$ with $4p=L^2+27M^2$ and
$L\e 1\mod 3$, $\sum_{n=0}^{p-1}\f{W_n}{(-3)^n}\e -L+\f pL\mod
{p^2}$. This was conjectured by the author's brother Z.W. Sun in
[Su1].

\par Throughout this paper, let $q_p(a)=(a^{p-1}-1)/p$
for given odd prime $p$ and $a\in\Bbb Z$, and let $H_0=0$ and
$H_n=1+\f 12+\cdots+\f 1n$ $(n\ge 1)$.

\section*{2. Congruences for $\sum_{k=0}^{p-1}\b
ak\b{-1-a}k\f 1{k+b}\mod {p^2}.$}
\par We begin with basic congruences for harmonic numbers $H_n$.
For odd prime $p$ and
$k\in\{1,2,\ldots,p-1\}$, it is clear that
$$H_{p-1-k}=\sum_{i=1}^{p-1}\f 1i-\sum_{i=1}^k\f 1{p-i}\e
0+\sum_{i=1}^k\f 1i=H_k\mod p.\tag 2.1$$ For $k=1,2,\ldots,\f{p-1}2$
we have
$$H_{\f{p-1}2-k}=H_{\f{p-1}2}-\sum_{i=1}^k\f 1{\f{p+1}2-i}
\e H_{\f{p-1}2}+2\sum_{i=1}^k\f 1{2i-1}=H_{\f{p-1}2}+2(H_{2k}-\f
12H_k)\mod p.$$ It is well known that (see [L])
$$H_{\f{p-1}2}\e -2q_p(2)\qtq{and}H_{[\f p3]}\e -\f 32q_p(3)\mod p
\ \t{for}\ p>3.\tag 2.2$$ Thus
$$2H_{2k}\e 2q_p(2)+H_k+H_{\f{p-1}2-k}\mod
p\qtq{for}k=1,2,\ldots,\f{p-1}2.\tag 2.3$$ From (2.2) and (2.3) we
deduce the known congruences (see [L]):
$$H_{[\f p4]}\e -3q_p(2)\mod p\qtq{and}H_{[\f p6]}
\e -2q_p(2)-\f 32q_p(3)\mod p\ \t{for}\ p>3.\tag 2.4$$
\par Let $p$ be an odd prime, $a,m\in\Bbb Z_p$, $k\in\{1,2,\ldots,p-1\}$ and $a-i\not\e 0\mod {p^2}$ for
$i=0,1,\ldots,k-1$. Then clearly
$$\align \b{a+mp}k&=\f{(mp+a)(mp+a-1)\cdots(mp+a-k+1)}{k!}
\\&\e \f
1{k!}\Big(a(a-1)\cdots(a-k+1)+mp\sum_{i=0}^{k-1}\f{a(a-1)
\cdots(a-k+1)}{a-i}\Big)\mod{p^2}.\endalign$$ That is,
$$\b{a+mp}k\e\b ak\Big(1+\sum_{i=0}^{k-1}\f{mp}{a-i}\Big)\mod
{p^2}.\tag 2.5$$

 \par By [S9,(2.2)],
for $n=0,1,2,\ldots$,
$$\sum_{k=0}^n\b
ak\b{-1-a}k\f {ab+1}{bk+1}-\sum_{k=0}^n\b {a-1}k\b{-a}k\f
{ab-1}{bk+1}=2\b{a-1}n\b{-a-1}n.$$ Replacing $b$ with $\f 1b$ yields
$$\sum_{k=0}^n\b ak\b{-1-a}k\f {a+b}{k+b}-\sum_{k=0}^n\b
{a-1}k\b{-a}k\f {a-b}{k+b}=2\b{a-1}n\b{-a-1}n.\tag 2.6$$ By [S9,
Lemma 2.2], for odd prime $p$ and $a\in\Bbb Z_p$ with $a\not\e 0\mod
p$,
$$\b {a-1}{p-1}\b{-a-1}{p-1}\e \f{(a-\ap)(p+a-\ap)}{a^2}\mod {p^3}.\tag 2.7$$
Now taking $b=a$ in (2.6) and applying (2.7) yields
$$\sum_{k=0}^{p-1}\f{\b ak\b{-1-a}k}{k+a}=\f
1a\b{a-1}{p-1}\b{-a-1}{p-1}\e \f{(a-\ap)(p+a-\ap)}{a^3}\mod
{p^3}.\tag 2.8$$
 \pro{Lemma 2.1} Let $p$ be an odd prime,
$a,b\in\Bbb Z_p$,  $a\not\e 0\mod p$, $b+\langle-b\rangle_p\not\e
0\mod {p^2}$ and
 $$S(a,b)=\sum_{k=0}^{p-1}\b
ak\b{-1-a}k\f p{k+b}.$$ Suppose that $m\in\{0,1,\ldots,\ap-1\}$ and
$\langle a+b\rangle_p\not\in\{0,1,\ldots,m\}$. Then
$$S(a,b)\e \f{\b{a-b}{m+1}}{\b{a+b}{m+1}}S(a-m-1,b)\mod{p^3}.$$
\endpro
Proof. Since $b+\langle-b\rangle_p\not\e 0\mod {p^2}$ we see that
$b+k\not\e 0\mod {p^2}$ and so $\f p{k+b}\in\Bbb Z_p$ for
$k=0,1,\ldots,p-1$. Hence $S(a,b)\in\Bbb Z_p$. For
$k=0,1,\ldots,\ap-1$ we have $a-k\not\e 0\mod p$. From (2.6) and
(2.7),
$$\align&(a-k+b)S(a-k,b)-(a-k-b)S(a-k-1,b)
\\&=2p\b{a-k-1}{p-1}\b{-(a-k)-1}{p-1}\e 0\mod {p^3}.\endalign$$
Thus, for $k=0,1,\ldots,m$,
$$S(a-k,b)\e \f{a-b-k}{a+b-k}S(a-k-1,b)\mod {p^3}.$$ Therefore,
$$\align S(a,b)&\e \f{a-b}{a+b}S(a-1,b)\e \f{a-b}{a+b}\cdot
\f{a-b-1}{a+b-1}S(a-2,b)\\&\e\cdots\e
\prod_{k=0}^m\f{a-b-k}{a+b-k}\cdot S(a-m-1,b)\mod{p^3}.\endalign$$
To see the result, we note that
$$\prod_{k=0}^m\f{a-b-k}{a+b-k}
=\f{\b{a-b}{m+1}}{\b{a+b}{m+1}}.$$

\pro{Lemma 2.2} Let $p$ be an odd prime, $b,t\in\Bbb Z_p$, $bt\not\e
0\mod p$ and $s=(b-\bp)/p\not\e -1\mod p$. Then
 $$\sum_{k=0}^{p-1}\b
{pt}k\b{-1-pt}k\f p{k+b}\e \Big(1+\f t{s+1}\Big)\f pb\Big(1-\f
{pt}b\Big) \mod {p^3}.$$
\endpro
Proof. It is well known that $\sum_{k=1}^{p-1}\f 1k\e 0\mod p$.
Thus,
$$\aligned &\sum_{k=1}^{p-1}\b
{pt}k\b{-1-pt}k\f p{k+b}\\&=\sum_{k=1}^{p-1}
(-1)^k\f{pt(pt+k)(p^2t^2-1^2)\cdots(p^2t^2-(k-1)^2)}{k!^2}\cdot\f
p{k+b}\\&\e -\sum_{k=1}^{p-1}\f{pt(pt+k)}{k^2}\cdot\f
p{k+b}=-p^2t^2\sum_{k=1}^{p-1}\f
p{k^2(k+b)}-\f{pt}b\sum_{k=1}^{p-1}\Big(\f pk-\f{p}{k+b}\Big)
\\&\e -\f{p^2t^2}{(p-\bp)^2(p-\bp+b)/p}+\f{pt}b\Big(\f 1{(p-\bp+b)/p}
+\sum\Sb k=1\\k\not=p-\bp\endSb^{p-1}\f p{k+b}\Big)
\\&\e -\f{p^2t^2}{b^2(s+1)}+\f{pt}b\Big(\f
1{s+1}-\f pb\Big)
 \mod{p^3},
\endaligned$$
which yields the result.
\par\q
 \pro{Theorem 2.1} Let $p$ be an odd prime,
$a,b\in\Bbb Z_p$,  $ab\not\e 0\mod p$ and $\bp\le p-1-\ap$. Assume
$t=(a-\ap)/p$ and $s=(b-\bp)/p\not\e -1\mod p$. For $\bp\le\ap$ we
have
$$\sum_{k=0}^{p-1}\b
ak\b{-1-a}k\f 1{k+b}\e
\f{p(s+t+1)(s-t)}{b^2(s+1)\b{\ap}{\bp}\b{p-1-\ap}{\bp}} \mod
{p^2}.$$ For $\bp>\ap$ we have
 $$\align\sum_{k=0}^{p-1}\b
ak\b{-1-a}k\f 1{k+b}&\e \f{s+1+t}{b(s+1)}\cdot
\f{\b{\bp-1}{\ap}}{\b{p-1-\bp}{\ap}}\Big(1
+p\f{s+1}b+p(2s+1)H_{\bp-1}\\&\q-p(s-t)H_{\bp-\ap-1}
-p(s+t+1)H_{\ap+\bp} \Big) \mod {p^2}.\endalign$$
\endpro
Proof. For $k\in\{0,1,\ldots,p-1\}$ clearly $k+b\not\e 0\mod p$ for
$k\not=p-\bp$. For $k=p-\bp$ we have $k+b\not\e 0\mod {p^2}$, $\b
ak=\b a{p-\bp}\e \b {\ap}{p-\bp}=0\mod p$ and so $\b ak\f
1{k+b}\in\Bbb Z_p$. Therefore, $\b ak\b{-1-a}k\f 1{k+b}\in\Bbb Z_p$
for $k=0,1,\ldots,p-1$. Now taking $m=\ap-1$ in Lemma 2.1 and then
applying Lemma 2.2 gives
$$\sum_{k=0}^{p-1}\b
ak\b{-1-a}k\f 1{k+b}\e \f{\b{a-b}{\ap}}{\b{a+b}{\ap}} \Big(1+\f
t{s+1}\Big)\f 1b\Big(1-\f {pt}b\Big)\mod {p^2}.\tag 2.9$$ We first
assume $\bp\le \ap$. Since $a-b\e \ap-\bp\mod p$, we see that
$$\align \f{\b{a-b}{\ap}}{\b{a+b}{\ap}}
&= (a-b-(\ap-\bp))\f{(a-b)(a-b-1)\cdots
(a-b-(\ap-\bp-1))}{\ap!}\\&\qq\times\f{(a-b-(\ap-\bp+1))
\cdots(a-b-\ap+1)}{\b{a+b}{\ap}}
\\&\e p(t-s)\f{(\ap-\bp)!}{\ap!}\cdot\f{(-1)(-2)\cdots(-(\bp-1))}{\b{\ap+\bp}{\ap}}
\\&=\f{p(t-s)}{-\bp}\cdot\f{(\ap-\bp)!(-1)^{\bp}\bp!}{\ap!(-1)^{\bp}\b{-1-\ap}{\bp}}
\e \f{p(s-t)}{b\b{\ap}{\bp}\b{p-1-\ap}{\bp}} \mod {p^2}.
\endalign$$ This together with (2.9) yields the result in the case $\bp\le
\ap$.
\par For $\bp>\ap$, using (2.5) we see that
$$\align \f{\b{a-b}{\ap}}{\b{a+b}{\ap}}
&=\f{\b{b-1-a+\ap}{\ap}}{\b{-b-1-a+\ap}{\ap}}=\f{\b{\bp-1+p(s-t)}{\ap}}
{\b{p-1-\bp-p(s+t+1)}{\ap}}
\\&
\e\f{\b{\bp-1}{\ap}(1+p(s-t)\sum_{i=0}^{\ap-1}\f 1{\bp-1-i})}
{\b{p-1-\bp}{\ap}(1-p(s+t+1)\sum_{i=0}^{\ap-1}\f 1{p-1-\bp-i})}
\\&\e\f{\b{\bp-1}{\ap}(1+p(s-t)(H_{\bp-1}-H_{\bp-\ap-1})}
{\b{p-1-\bp}{\ap}(1+p(s+t+1)(H_{\ap+\bp}-H_{\bp})}
\\&\e\f{\b{\bp-1}{\ap}}{\b{p-1-\bp}{\ap}}
(1+p(s-t)(H_{\bp-1}-H_{\bp-\ap-1})\\&\qq\times(1-p(s+t+1)(H_{\ap+\bp}-H_{\bp})
\\&\e\f{\b{\bp-1}{\ap}}{\b{p-1-\bp}{\ap}}
\Big(1+p\f{s+t+1}b+p(2s+1)H_{\bp-1}\\&\qq-p(s-t)H_{\bp-\ap-1}
-p(s+t+1)H_{\ap+\bp}\Big)
 \mod {p^2}.
\endalign$$
This together with (2.9) yields the remaining part.
\par\q
 \pro{Corollary 2.1} Let $p$ be a prime with $p>3$, $a\in\Bbb
Z_p$, $1\le \ap\le \f{p-3}2$ and $t=(a-\ap)/p\not\e -1\mod p$. Then
$$\align &\sum_{k=0}^{p-1}\b ak\b{-1-a}k\f 1{k+a+1}
\\&\e\f{2t+1}{(a+1)(t+1)\b {p-2-\ap}{\ap}}
\Big(1+p\f{t+1}{a+1}-p(2t+1)(H_{2\ap+1}-H_{\ap})\Big)
 \mod {p^2}.\endalign$$
 \endpro
 Proof. Set $b=a+1$. Then $\bp=\ap+1$ and so $\ap<\bp\le p-1-\ap$.
 Now putting $b=a+1$ in Theorem 2.1 and noting that $s=t$ yields the result.
\par\q
\pro{Theorem 2.2} Let $p>3$ be a prime. Suppose $b\in\Bbb Z_p$,
 $\bp\not=0$ and $s=(b-\bp)/p\not\e -1\mod p$. Then
 $$\align &\sum_{k=0}^{p-1}\f{\b{2k}k^2}{16^k(k+b)}
 \e \f {(s+1/2)^2p}{b^2(s+1)\b{(p-1)/2}{\bp}^2}
\mod{p^2}\qtq{for}\bp<\f p2,\tag 2.10
\\&\sum_{k=0}^{(p-1)/2}\f{\b{2k}k^2}{16^k(k+b)}
 \e \f {b-\bp}{b^2\b{(p-1)/2}{\bp}^2}
\mod{p^2}\qtq{for}\bp<\f p2,\tag 2.11
\\&\sum_{k=0}^{p-1}\f{\b{2k}k\b{3k}k}{27^k(k+b)}\e (-1)^{\bp}\f
{(s+1/3)(s+2/3)p}{b^2(s+1)\b{2\bp}{\bp}\b{[p/3]+\bp}{[p/3]-\bp}}
\mod{p^2}\ \t{for}\ \bp<\f p3,\tag 2.12
\\&\sum_{k=0}^{p-1}\f{\b{2k}k\b{4k}{2k}}{64^k(k+b)}\e (-1)^{\bp}\f
{(s+1/4)(s+3/4)p}{b^2(s+1)\b{2\bp}{\bp}\b{[p/4]+\bp}{[p/4]-\bp}}
\mod{p^2}\ \t{for}\ \bp<\f p4,\tag 2.13
\\&\sum_{k=0}^{p-1}\f{\b{3k}k\b{6k}{3k}}{432^k(k+b)}\e (-1)^{\bp}\f
{(s+1/6)(s+5/6)p}{b^2(s+1)\b{2\bp}{\bp}\b{[p/6]+\bp}{[p/6]-\bp}}
\mod{p^2}\ \t{for}\ \bp<\f p6.\tag 2.14
\endalign$$
\endpro
Proof. Set $a=-\f 12$. Then $\ap=\f{p-1}2$ and $t=(a-\ap)/p=-1/2$.
For $\bp<\f p2$ we have $\bp\le p-1-\ap=\ap$. Since
$\b{-1/2}k^2=\b{2k}k16^{-k}$, taking $a=-\f 12$ in Theorem 2.1
yields (2.10).  Clearly,
$$\align\f 1p\b{2(p-\bp)}{p-\bp}
&=\f{(2p-2\bp)(2p-(2\bp-1))\cdots(p+1)\cdot(p-1)\cdots(p-(\bp-1))}{(p-\bp)!}
\\&\e \f{(p-2\bp)!(-1)^{\bp-1}(\bp-1)!}{(p-\bp)!}
\e \f{2(-1)^{\bp}\bp!}{(p-\bp))\cdots (p-2\bp)}
\\&\e
\f{2}{-b\b{2\bp}{\bp}}\e \f{2}{-b(-4)^{\bp}\b{(p-1)/2}{\bp}}
\mod{p}.\endalign$$ Since $p\mid \b{2k}k$ for $\f p2<k<p$, using
(2.10) we see that
$$\align\sum_{k=0}^{(p-1)/2}\f{\b{2k}k^2}{16^k(k+b)}
&\e\sum_{k=0}^{p-1}\f{\b{2k}k^2}{16^k(k+b)}- \f
{\b{2(p-\bp)}{p-\bp}^2}{16^{p-\bp}(p-\bp+b)}
\\&\e\f {p(s+1/2)^2}{b^2(s+1)\b{(p-1)/2}{\bp}^2}-
\f 1{16^{p-\bp}(s+1)}\cdot\f {4p}{16^{\bp}b^2\b{(p-1)/2}{\bp}^2}
\\&\e\f{ps}{b^2\b{(p-1)/2}{\bp}^2}=\f {b-\bp}{b^2\b{(p-1)/2}{\bp}^2}
\mod {p^2}.\endalign$$ This proves (2.11).
\par Now we consider (2.12)-(2.14). We first note that
$$\align \b{\ap}{\bp}\b{p-1-\ap}{\bp}&\e \b{\ap}{\bp}\b{-1-\ap}{\bp}
=(-1)^{\bp}\b{\ap}{\bp}\b{\ap+\bp}{\ap}\\&=(-1)^{\bp}\b{2\bp}{\bp}
\b{\ap+\bp}{\ap-\bp} \mod p.\endalign$$
 Set $a=-\f 13$ or $-\f 23$ according as
 $p\e 1\mod 3$ or $p\e 2\mod 3$. Then  $\ap=[\f p3]$ and $t=(a-\ap)/p=
 -\f 13$. Recall that
 $\b{-1/3}k\b{-2/3}k=\b{2k}k\b{3k}k27^{-k}$. From Theorem 2.1 and the above
 we obtain (2.12).
\par  Next consider (2.13). Set $a=-\f 14$ or $-\f 34$ according as
 $p\e 1\mod 4$ or $p\e 3\mod 4$. Then  $\ap=[\f p4]$ and $t=(a-\ap)/p=
 -\f 14$. Recall that
 $\b{-1/4}k\b{-3/4}k=\b{2k}k\b{4k}{2k}64^{-k}$, (2.13) follows from
 Theorem 2.1 (with $a=-\f 14$).
\par Finally consider (2.14). Set $a=-\f 16$ or $-\f 56$ according as
 $p\e 1\mod 6$ or $p\e 5\mod 6$. Then  $\ap=[\f p6]$ and $t=(a-\ap)/p=
 -\f 16$. Since
 $\b{-1/6}k\b{-5/6}k=\b{3k}k\b{6k}{3k}432^{-k}$, applying Theorem 2.1
  yields (2.14). The proof is now complete.
  \par\q
  \par Let $p$ be a prime with $p>3$. As examples,
  from (2.10), (2.11) and the fact $\b{(p-1)/2}k\e
  \b{-1/2}k=\b{2k}k4^{-k}\mod p$   we deduce
  that
  $$\align &-\sum_{k=0}^{(p-1)/2}\f{\b{2k}k^2}{16^k(4k+3)}
\e 3\sum_{k=0}^{p-1}\f{\b{2k}k^2}{16^k(4k+3)} \e
p\b{\f{p-1}2}{\f{p-1}4}^{-2}\mod{p^2}\ \t{for}\ p\e 1 \mod 4,
\\&8\sum_{k=0}^{p-1}\f{\b{2k}k^2}{16^k(3k+1)} \e
5\sum_{k=0}^{p-1}\f{\b{2k}k^2}{16^k(6k+1)} \e
p\b{\f{p-1}2}{\f{p-5}6}^{-2}\mod {p^2}\ \t{for}\ p\e 5\mod
6,\\&45\sum_{k=0}^{p-1}\f{\b{2k}k^2}{16^k(8k+1)} \e
7\sum_{k=0}^{p-1}\f{\b{2k}k^2}{16^k(8k+3)}\e -p\b{[\f p4]}{[\f
p8]}^{-2}\mod {p^2}\ \t{for}\ p\e 5\mod 8,
\\&\f{14}9\sum_{k=0}^{p-1}\f{\b{2k}k^2}{16^k(8k+1)} \e
10\sum_{k=0}^{p-1}\f{\b{2k}k^2}{16^k(8k+3)}\e p\b{[\f p4]}{[\f
p8]}^{-2}\mod {p^2}\ \t{for}\ p\e 7\mod 8.
\endalign$$

\pro{Theorem 2.3} Let $p$ be a prime with $p>3$. Then
$$\align &\sum_{k=0}^{p-1}\f{\b{2k}k\b{3k}k}{27^k(3k+2)}
\e \f 1{2\b{2(p-1)/3}{(p-1)/3}}\mod {p^2}\qtq{for}p\e 1\mod 3,\tag
2.15 \\&\sum_{k=0}^{p-1}\f{\b{2k}k\b{3k}k}{27^k(3k+1)} \e \f
{1+2p}{2\b{2(p-2)/3}{(p-2)/3}}\mod {p^2}\qtq{for}p\e 2\mod 3,\tag
2.16
\\&\sum_{k=0}^{p-1}\f{\b{2k}k\b{4k}{2k}}{64^k(4k+3)}\e
(-1)^{\f{p-1}4}\f{2^{p-1}+1}{6\b{(p-1)/2}{(p-1)/4}}\mod {p^2} \q
\t{for}\q p\e 1\mod 4,\tag 2.17
\\&\sum_{k=0}^{p-1}\f{\b{2k}k\b{4k}{2k}}{64^k(4k+1)}\e
(-1)^{\f{p+1}4}\f{2^{p-1}+2p+1}{6\b{(p-3)/2}{(p-3)/4}}\mod {p^2} \q
\t{for}\q p\e 3\mod 4,\tag 2.18
\\&\sum_{k=0}^{p-1}\f{\b{3k}k\b{6k}{3k}}{432^k(6k+5)}\e
\f{1-\f 43(2^{p-1}-1)}{5\b{5(p-1)/6}{(p-1)/6}}\mod {p^2} \q
\t{for}\q p\e 1\mod 6,\tag 2.19
\\&\sum_{k=0}^{p-1}\f{\b{3k}k\b{6k}{3k}}{432^k(6k+1)}\e
\f 45\cdot\f{1+5p-\f 43(2^{p-1}-1)}{\b{(5p-7)/6}{(p-5)/6}}\mod {p^2}
\ \t{for}\ p\e 5\mod 6.\tag 2.20
\endalign$$
\endpro
Proof.  For $a=-\f 13$ and $p\e 1\mod 3$ we see that $\ap=\f{p-1}3$,
$\f{a-\ap}p=-\f 13$ and
  $$H_{2\ap+1}-H_{\ap}=H_{\f{2p+1}3}-H_{\f{p-1}3}
  \e H_{p-1-\f{2p+1}3}-H_{\f{p-1}3}=H_{\f{p-4}3}
  -H_{\f{p-1}3}=-\f 3{p-1}\e 3\mod p.$$ Thus, from Corollary
  2.1 and (1.2),
 $$\align \sum_{k=0}^{p-1}\f{\b{2k}k\b{3k}k}{27^k(3k+2)}
 &=\f 13\sum_{k=0}^{p-1}\b{-\f 13}k\b{-\f 23}k\f 1{k+2/3}
\\&\e\f{2t+1}{3(a+1)(t+1)\b {p-2-(p-1)/3}{(p-1)/3}}
\e \f{1}{2\b{2(p-1)/3}{(p-1)/3}}\mod{p^2},
\endalign$$
For $a=-\f 23$ and $p\e 2\mod 3$ we see that $\ap=\f{p-2}3$,
$\f{a-\ap}p=-\f 13$ and
$$H_{2\ap+1}-H_{\ap}\e
H_{p-2-2\ap}-H_{\ap}=H_{\f{p-2}3}-H_{\f{p-2}3} =0\mod p.$$ Thus,
from Corollary 2.1 and (1.2),
 $$\align \sum_{k=0}^{p-1}\f{\b{2k}k\b{3k}k}{27^k(3k+1)}
 &=\f 13\sum_{k=0}^{p-1}\b{-\f 13}k\b{-\f 23}k\f 1{k+1/3}
\\&\e\f{2t+1}{3(a+1)(t+1)\b {p-2-(p-2)/3}{(p-2)/3}}
\Big(1+p\f{t+1}{a+1}\Big) =\f
{1+2p}{2\b{2(p-2)/3}{(p-2)/3}}\mod{p^2}.
\endalign$$
\par  For $a=-\f 14$ and $p\e 1\mod 4$
 we see that $\ap=\f{p-1}4$ and
$\f{a-\ap}p=-\f 14$. It is well known that $H_{\f{p-1}2}\e
-2q_p(2)\mod p$ and $H_{\f{p-1}4}\e -3q_p(2)\mod p$. Thus,
$$H_{2\ap+1}-H_{\ap}=\f 1{(p+1)/2}+H_{\f{p-1}2}-H_{\f{p-1}4}\e
2-2q_p(2)+3q_p(2)=2+q_p(2)\mod p.$$ By [S3, Lemma 2.5],
$$(-1)^{\f{p-1}4}\b{3(p-1)/4}{(p-1)/4}
\e \Big(3-2(-4)^{\f{p-1}4}\Big)\b{(p-1)/2}{(p-1)/4} \e (1-pq_p(2))
\b{(p-1)/2}{(p-1)/4} \mod {p^2}.$$ Therefore, taking $a=-\f 14$ in
Corollary 2.1 and applying the above and (1.2) gives
$$\align&\sum_{k=0}^{p-1}\f{\b{2k}k\b{4k}{2k}}{64^k(4k+3)}
\\&=\f 14\sum_{k=0}^{p-1}\b{-\f 14}k\b{-\f 34}k\f 1{k+3/4} \e\f
14\cdot\f {\f 12}{\f 34\cdot\f 34\b{p-2-(p-1)/4}{(p-1)/4}}
\Big(1+p-\f p2(H_{\f{p-1}2+1}-H_{\f{p-1}4})\Big)
\\&\e \f 1{3\b{3(p-1)/4}{(p-1)/4}}\Big(1-\f p2q_p(2)\Big)
\e \f 13(-1)^{\f{p-1}4}\f 1{(1-pq_p(2))\b{(p-1)/2}{(p-1)/4}}
\Big(1-\f p2q_p(2)\Big)
\\&\e \f 13(-1)^{\f{p-1}4}\f
1{\b{(p-1)/2}{(p-1)/4}}(1+pq_p(2))\Big(1-\f p2q_p(2)\Big) \e  \f
13(-1)^{\f{p-1}4}\f {1+pq_p(2)/2}{\b{(p-1)/2}{(p-1)/4}}
 \mod
{p^2}.\endalign$$
\par  For $a=-\f 34$ and $p\e 3\mod 4$
 we see that $\ap=\f{p-3}4$ and
$\f{a-\ap}p=-\f 14$. By [S3, Lemma 2.5],
$$\align (-1)^{\f{p-3}4}\f{\b{(p-1)/2+(p-3)/4}{(p-3)/4}}
{\b{(p-1)/2}{(p-3)/4}} &\e
3-2(-4)^{\f{p-3}4}\f{\b{(p-1)/2}{(p-3)/4}}
{\b{(p-3)/2}{(p-3)/4}}=3-2(-4)^{\f{p-3}4}\f{2p-2}{p+1} \\& \e
3-2\cdot(-1)^{\f{p+1}4}2^{\f{p-1}2}(1-2p) \\&\e 3-2(1+pq_p(2)/2)+4p
=1+4p-pq_p(2)\mod{p^2}.\endalign$$ Hence
$$\align &(-1)^{\f{p-3}4}\b{\f{p-1}2+\f{p-3}4}{\f{p-3}4}\e
(1+4p-pq_p(2)){\b{\f{p-1}2}{\f{p-3}4}}
\\&=(1+4p-pq_p(2))\f{2p-2}{p+1}{\b{(p-3)/2}{(p-3)/4}}
\e -2(1+2p-pq_p(2)){\b{(p-3)/2}{(p-3)/4}}
\endalign $$
Also, $$H_{2\ap+1}-H_{\ap}=H_{\f{p-1}2}-H_{\f{p-3}4}
\e-2q_p(2)+3q_p(2)=q_p(2)\mod p.$$ Now, taking $a=-\f 34$ in
Corollary 2.1 and then applying the above gives
$$\align\sum_{k=0}^{p-1}\f{\b{2k}k\b{4k}{2k}}{64^k(4k+1)}
&=\f 14\sum_{k=0}^{p-1}\b{-\f 14}k\b{-\f 34}k\f 1{k+1/4} \e\f
23\cdot\f{1+3p-pq_p(2)/2}{\b{\f{p-1}2+\f{p-3}4}{\f{p-3}4}}
\\&\e \f 23(-1)^{\f{p+1}4}\f{1+3p-\f 12pq_p(2)}{2(1+2p-pq_p(2))
\b{\f{p-3}2}{\f{p-3}4}} \e \f 13(-1)^{\f{p+1}4}\f{1+p+\f
12pq_p(2)}{\b{\f{p-3}2}{\f{p-3}4}}\mod {p^2}.
\endalign$$
\par For $a=-\f 16$ and $p\e 1\mod 6$ we have $\ap=\f{p-1}6$,
$\f{a-\ap}p =-\f 16$. Applying (2.2) and (2.4),
$$H_{2\ap+1}-H_{\ap}=\f 1{(p+2)/3}+H_{\f{p-1}3}-H_{\f{p-1}6}
\e\f 32+2q_p(2)\mod p.$$ Now taking $a=-\f 16$ in Corollary 2.1
yields
$$\align&\sum_{k=0}^{p-1}\f{\b{3k}k\b{6k}{3k}}{64^k(6k+5)}
=\f 16\sum_{k=0}^{p-1}\b{-\f 16}k\b{-\f 56}k\f 1{k+5/6}\\& \e \f{\f
16\cdot\f 23}{\sls 56^2\b{p-2-(p-1)/6}{(p-1)/6}}\Big(1+p-\f
23p\Big(\f 32+2q_p(2)\Big)\Big)\e\f {1-\f
43(2^{p-1}-1)}{5\b{5(p-1)/6}{(p-1)/6}} \mod{p^2}.\endalign$$
\par For $a=-\f 56$ and $p\e 5\mod 6$ we have $\ap=\f{p-5}6$ and $
\f{a-\ap}p=-\f 16$. Also, $H_{2\ap+1}-H_{\ap}=H_{[\f p3]}-H_{[\f
p6]}\e 2q_p(2)\mod p.$ Now taking $a=-\f 56$ in Corollary 2.1 and
applying the above yields (2.20). The proof is now complete.
\par\q
\pro{Theorem 2.4} Let $p$ be a prime with $p>3$. Then
$$\align&\sum_{k=0}^{p-1}
\f{\b{3k}k\b{6k}{3k}}{432^k(3k+1)}\e \cases\f{2^{p-1}+2}{3\cdot
2^{\f{p-1}3}}\mod {p^2}&\t{if $p\e 1\mod 6$,}
\\2^{\f{p+1}3}(2^{p-3}-1)\mod {p^2}&\t{if $p\e 5\mod 6$.}
\endcases
\endalign$$
\endpro
Proof.  We first assume $p\e 1\mod {6}$. Set $a=-\f 16$ and $b=\f
13$. Then $\ap=\f{p-1}6$, $\bp=\f{2p+1}3$, $t=(a-\ap)/p=-\f 16$ and
$s=(b-\bp)/p=-\f 23$. Since $\ap<\bp<p-1-\ap$, from Theorem 2.1 and
(2.1)-(2.4),
$$\align \sum_{k=0}^{p-1}
\f{\b{3k}k\b{6k}{3k}}{432^k(3k+1)}&\e  \f{\b{(2p-2)/3}{(p-1)/6}}
{2\b{(p-4)/3}{(p-1)/6}}\Big(1+p-\f p3H_{\f{2p-2}3}+\f
12pH_{\f{p-1}2} -\f p6H_{\f{5p+1}{6}}\Big)
\\&\e
\f{\b{(2p-2)/3}{(p-1)/6}} {2\b{(p-4)/3}{(p-1)/6}}\Big(1+p-\f
p3H_{\f{p-1}3}+\f 12pH_{\f{p-1}2} -\f p6\big(H_{\f{p-1}{6}}-\f
1{(p-1)/6}\big)\Big)
\\&\e
\f{\b{(2p-2)/3}{(p-1)/6}} {\b{(p-1)/3}{(p-1)/6}}\Big(1-\f
23pq_p(2)+\f 34pq_p(3)\Big) \mod {p^2}\endalign$$ By [S3, Lemmas 2.4
and 2.5], for $k=1,2,\ldots,\f{p-1}2$,
$$\align\b{\f{p-1}2+k}k&\e
(-1)^k\b{\f{p-1}2}k\Big(1+2p\sum_{i=1}^k\f 1{2i-1}\Big)
\\&\e \b{2k}k4^{-k}\Big(1-p\sum_{i=1}^k\f 1{2i-1}\Big)\Big(1+2p\sum_{i=1}^k\f 1{2i-1}\Big)
\\&\e \b{2k}k 4^{-k}\Big(1+p\sum_{i=1}^k\f 1{2i-1}\Big)
\mod {p^2}. \endalign$$ That is,
$$\b{\f{p-1}2+k}k\e \b{2k}k
4^{-k}\Big(1+p\big(H_{2k}-\f 12H_k\big)\Big)\mod {p^2}.\tag 2.21$$
Thus,
$$\f{\b{(2p-2)/3}{(p-1)/6}}{\b{(p-1)/3}{(p-1)/6}}\e
4^{-\f{p-1}6}\Big(1+p\big(H_{\f{p-1}3}-\f 12H_{\f{p-1}6}\big)\Big)
\e 2^{-\f{p-1}3}\Big(1+pq_p(2)-\f 34pq_p(3)\Big)\mod {p^2}.$$
Therefore,
$$\align \sum_{k=0}^{p-1}
\f{\b{3k}k\b{6k}{3k}}{432^k(3k+1)}&\e \Big(1-\f 23pq_p(2)+\f
34pq_p(3)\Big)2^{-\f{p-1}3}\Big(1+pq_p(2)-\f 34pq_p(3)\Big)
\\&\e 2^{-\f{p-1}3}\Big(1+\f 13pq_p(2)\Big)=\f{2^{p-1}+2}{3\cdot
2^{\f{p-1}3}}
 \mod {p^2}.\endalign$$
 \par Now  assume $p\e 5\mod {6}$. Set $a=-\f 56$ and $b=\f
13$. Then $\ap=\f{p-5}6$, $\bp=\f{p+1}3$, $t=(a-\ap)/p=-\f 16$ and
$s=(b-\bp)/p=-\f 13$. Since $\ap<\bp<p-1-\ap$, from Theorem 2.1 and
(2.1)-(2.4),
$$\align &\sum_{k=0}^{p-1}
\f{\b{3k}k\b{6k}{3k}}{432^k(3k+1)}\\&\e  \f
34\cdot\f{\b{(p-2)/3}{(p-5)/6}} {\b{(2p-4)/3}{(p-5)/6}}\Big(1+2p+\f
p3H_{\f{p-2}3}+\f p6H_{\f{p+1}6} -\f p2H_{\f{p-1}{2}}\Big)
\\&\e  \f
34\cdot\f{\b{(p-2)/3}{(p-5)/6}} {\b{(2p-4)/3}{(p-5)/6}}\Big(1+2p+\f
p3\Big(-\f 32\qp 3\Big)+\f p6\Big(\f 1{(p+1)/6}-2\qp 2-\f 32\qp
3\Big) -\f p2(-2\qp 2)\Big)
\\&\e\f 34\cdot \f{2(p-2)}{p+1}\cdot\f{\b{(p-5)/3}{(p-5)/6}}
 {\b{(2p-4)/3}{(p-5)/6}}\Big(1+3p+\f 23p\qp 2-\f 34p\qp 3\Big)
\mod {p^2}.\endalign$$ By (2.21),
$$\align \b{\f{2p-4}3}{\f{p-5}6}&
\e \b{\f{p-5}3}{\f{p-5}6}4^{-\f{p-5}6}\Big(1+p\Big(H_{\f{p-5}3}-\f
12H_{\f{p-5}6}\Big)\Big)
\\&\e \b{\f{p-5}3}{\f{p-5}6}4^{-\f{p-5}6}\Big(1+p\Big(-\f
32\qp 3-\f 1{(p-2)/3}-\f 12\Big(-2\qp 2-\f 32\qp 3\Big)\Big)\Big)
\\&\e  \b{\f{p-5}3}{\f{p-5}6}2^{-\f{p-5}3}\Big(1+p\Big(\f
32+\qp 2-\f 34\qp 3\Big)\Big)\mod {p^2}.\endalign$$
Hence
$$\align &\sum_{k=0}^{p-1}
\f{\b{3k}k\b{6k}{3k}}{432^k(3k+1)}\\&\e \f 32\cdot
\f{(p-2)(1-p)}{1-p^2}\cdot\f{2^{(p-5)/3}}{1+p(\f 32+\qp 2-\f 34\qp
3)}\Big(1+p\Big(3+\f 23\qp 2-\f 34\qp 3\Big)\Big)
\\&\e \f 32(-2+3p)2^{\f{p-5}3}\Big(1-p\Big(\f 32+\qp 2-\f 34\qp
3\Big)\Big)\Big(1+p\Big(3+\f 23\qp 2-\f 34\qp 3\Big)\Big)
\\&\e \f 32\cdot 2^{\f{p-5}3}\Big(-2+\f 23p\qp 2\Big)
=2^{\f{p+1}3}(2^{p-3}-1)\mod {p^2}.\endalign$$ This completes the
proof.

\pro{Conjecture 2.1} Let $p$ be a prime with $p\e 1\mod 4$. Then
$$\sum_{k=0}^{p-1} \f{\b{3k}k\b{6k}{3k}}{432^k(4k+1)}\e\cases
(-1)^y \mod {p^2}&\t{if $p=x^2+9y^2\e
1\mod{12}$,}\\(-3)^{\f{p-1}4}\mod {p}&\t{if $p\e
5\mod{12}$.}\endcases$$
\endpro
\pro{Conjecture 2.2} Let $p$ be a prime with $p>3$. Then
$$\sum_{k=0}^{p-1}\f{\b{2k}k\b{4k}{2k}}{64^k(3k+1)}
\e\cases 1\mod{p^2}&\t{if $p\e 1\mod 3$,}
\\-\f 54\mod {p^2}&\t{if $p\e 2\mod 3$}
\endcases$$
and
$$\sum_{k=0}^{p-1}\f{\b{2k}k\b{4k}{2k}}{64^k(3k+2)}
\e\cases \f 12\mod{p^2}&\t{if $p\e 1\mod 3$,}
\\-\f 25\mod {p^2}&\t{if $p\e 2\mod 3$.}
\endcases$$
\endpro
\pro{Conjecture 2.3} Let $p$ be a prime with $p\e 1\mod 6$ and so
$p=x^2+3y^2$. Then
$$\sum_{k=0}^{p-1}\f{\b{2k}k\b{4k}{2k}}{64^k}\cdot\f p{6k+1}
\e (-1)^{\f{p-1}6}\Big(4x^2-2p-\f{p^2}{4x^2}\Big)\mod{p^3}.$$
\endpro

\section*{3. The congruence for $\sum_{k=0}^{p-1}\b ak\b{-1-a}k\f p{k-a}
\mod {p^2}$}
\par Let $p$ be an odd prime. Replacing $a$ with $-1-a$ in Corollary
2.1 we obtain a congruence for $\sum_{k=0}^{p-1}\b ak\b{-1-a}k\f
1{k-a}\mod {p^2}$ under the condition $\f p2<\ap<p-1$. Now we
present the congruence for $\sum_{k=0}^{p-1}\b ak\b{-1-a}k\f
p{k-a}\mod {p^2}$ under the condition $1\le \ap<\f p2$. We note that
for $k=\ap$, $\b ak\b{-1-a}k\f 1{k-a}\not\in\Bbb Z_p$, but $\b
ak\b{-1-a}k\f p{k-a}\in\Bbb Z_p$ provided that $a-\ap\not\e
0\mod{p^2}$.

 \pro{Theorem 3.1} Let $p$ be an odd prime, $a\in\Bbb
Z_p$, $1\le \ap<\f p2$ and $t=\f{a-\ap}p\not\e 0\mod p$. Then
$$\align\sum_{k=0}^{p-1}\b ak\b{-1-a}k\f p{k-a}
&\e\f {2t+1}t(-1)^{\ap-1}\b{2\ap}{\ap}+2\b{p-1-\ap}{\ap}
 \\&\e (-1)^{\ap-1}\b{2\ap}{\ap}\Big(\f
1t+2p(H_{2\ap}-H_{\ap})\Big)\mod {p^2}.\endalign$$
\endpro
Proof. Set $S(a)=\sum_{k=0}^{p-1}\b ak\b{-1-a}k\f p{k-a}.$
 From [G, (1.41)],
$$\sum_{k=0}^n\b nk(-1)^k\f 1{k+x}=\f 1{x\b{x+n}n}=\f{(-1)^n}{(x+n)\b{-x}n}.
\tag 3.1$$ Thus, $\sum_{r=0}^k\b kr(-1)^r\f p{r-a} = \f
{(-1)^kp}{(k-a)\b ak}.$ From [S7, Theorem 2.4] and the above,
$$\aligned &(-1)^{\ap}S(a)\\&\e\sum_{k=0}^{p-1}\b ak\b{-1-a}k
\f {(-1)^kp}{(k-a)\b ak}=\sum_{k=0}^{p-1}(-1)^k\b{-1-a}k\f p{k-a}
\\&=\sum\Sb k=0\\k\not=\ap\endSb^{p-1}(-1)^k\b{-1-a}k\f p{k-a}
+(-1)^{\ap}\b{-1-a}{\ap}\f 1{(\ap-a)/p}
\\&\e\sum\Sb k=0\\k\not=\ap\endSb^{p-1-\ap}(-1)^k\b{p-1-\ap}k\f p{k-a}
+(-1)^{\ap}\b{-1-a}{\ap}\f 1{(\ap-a)/p}
\\&=\sum_{k=0}^{p-1-\ap}(-1)^k\b{p-1-\ap}k\f
p{k-a}+(-1)^{\ap}\b{p-1-\ap}{\ap}\f 1t -(-1)^{\ap}\b{-1-a}{\ap}\f 1t
\\&=\f{(-1)^{p-1-\ap}p}{(p-1-\ap-a)\b a{p-1-\ap}}
+\f{(-1)^{\ap}}t\Big(\b{p-1-\ap}{\ap}-\b{-1-a}{\ap}\Big) \mod{p^2}.
\endaligned$$
 Applying (2.5) we see that
$$\align S(a)&\e-\f p{(a+\ap-p+1)\b a{p-1-\ap}}
+\f 1t\Big(\b{-1-a+p(t+1)}{\ap}-\b{-1-a}{\ap}\Big)
\\&\e-\f p{(a+\ap-p+1)\b a{p-1-\ap}}
+\f 1t\b{-1-a}{\ap}\sum_{i=0}^{\ap-1}\f{p(t+1)}{-1-a-i}
\\&=-\f{(p-1-\ap)!\cdot p}{a(a-1)\cdots(a-(p-1-\ap))}-
\f {t+1}t\b{-1-a}{\ap}\sum_{i=0}^{\ap-1}\f p{1+a+i}
\mod{p^2}.\endalign$$
Since $\ap<\f p2$,
$$\align &\f{(p-1-\ap)!\cdot p}{a(a-1)\cdots(a-(p-1-\ap))}
\\&=\f{(p-1-\ap)!\cdot p}{(pt+\ap)(pt+\ap-1)
\cdots(pt+1)pt(pt-1)\cdots(pt-(p-1-2\ap))}
\\&\e\f {(p-1-\ap)!}{t\cdot \ap!(1+pt\sum_{i=1}^{\ap}\f 1i)
\cdot (p-1-2\ap)!(1-pt\sum_{i=1}^{p-1-2\ap}\f 1i)}
\\&\e \f 1t\b{p-1-\ap}{\ap}\Big(1-pt\sum_{i=1}^{\ap}\f 1i\Big)
\Big(1+pt\sum_{i=1}^{p-1-2\ap}\f 1i\Big)
\\&\e \f 1t\b{p-1-\ap}{\ap}\Big(1-ptH_{\ap}
+ptH_{p-1-2\ap}\Big)\mod {p^2}
\endalign$$
and
$$\align &\b{-1-a}{\ap}\sum_{i=0}^{\ap-1}\f p{1+a+i}
\\&\e \b{-1-a}{\ap}p\sum_{i=0}^{\ap-1}\f 1{\ap+1+i} \e
p\b{p-1-\ap}{\ap}(H_{2\ap}-H_{\ap})\mod{p^2}.\endalign$$ Recall that
$H_{p-1-k}\e H_k\mod p$. From the above we deduce that
$$\align &S(a)\\&\e -\f 1t\b{p-1-\ap}{\ap}(1+pt(H_{2\ap}-H_{\ap}))
-\f{p(t+1)}t\b{p-1-\ap}{\ap}(H_{2\ap}-H_{\ap})
\\&=-\f 1t\b{p-1-\ap}{\ap}-p\f {2t+1}t\b{p-1-\ap}{\ap}(H_{2\ap}-H_{\ap})\mod {p^2}.
\endalign$$
By (2.5),
$$ \b{p-1-\ap}{\ap}=(-1)^{\ap}\b{2\ap-p}{\ap}
\e(-1)^{\ap}\b{2\ap}{\ap}(1-p(H_{2\ap}-H_{\ap}))\mod {p^2}.$$ Hence
$$\align &S(a)\\&\e -\f
1t(-1)^{\ap}\b{2\ap}{\ap}(1-p(H_{2\ap}-H_{\ap})) -p\f
{2t+1}t(-1)^{\ap}\b{2\ap}{\ap}(H_{2\ap}-H_{\ap}) \\&=-\f 1t
(-1)^{\ap}\b{2\ap}{\ap}-2p(-1)^{\ap}\b{2\ap}{\ap}(H_{2\ap}-H_{\ap})
\\&\e \f 1t
(-1)^{\ap-1}\b{2\ap}{\ap}+2\Big(\b{p-1-\ap}{\ap}-(-1)^{\ap}\b{2\ap}{\ap}\Big)
 \mod {p^2}.\endalign$$ This completes the
proof.
\par\q
 \pro{Theorem 3.2} Let $p$ be a prime with $p\e 1\mod 4$. Then
$$\sum_{k=0}^{p-1}\f{\b{2k}k\b{4k}{2k}}{64^k}\cdot\f p{4k+1}
\e (-1)^{\f{p-1}4}\f{3-2^{p-1}}2\b{\f{p-1}2}{\f{p-1}4}\mod {p^2}.$$
\endpro
Proof. By (2.2) and (2.4), $H_{\f{p-1}2}\e -2q_p(2)\mod p$ and
$H_{[\f p4]}\e -3q_p(2)\mod p.$  Taking $a=-\f 14$ in Theorem 3.1
and noting that $\ap=\f{p-1}4$ and $t=-\f 14$ gives
$$\align&\sum_{k=0}^{p-1}\f{\b{2k}k\b{4k}{2k}}{64^k}\cdot\f p{4k+1}
=\f 14\sum_{k=0}^{p-1}\b{-\f 14}k\b{-\f 34}k\f p{k-(-\f 14)}
\\&\e (-1)^{\f{p-1}4}
\b{\f{p-1}2}{\f{p-1}4}(1-\f p2(H_{\f{p-1}2} -H_{\f{p-1}4})) \e
(-1)^{\f{p-1}4} \b{\f{p-1}2}{\f{p-1}4}(1-\f p2q_p(2))\mod {p^2}
\endalign$$
This proves the theorem.
\par\q
\pro{Theorem 3.3} Let $p$ be a prime with $p\e 1\mod 3$. Then
$$\sum_{k=0}^{p-1}\f{\b{2k}k\b{3k}k}{27^k}\cdot\f p{3k+1}\e
\b{2(p-1)/3}{(p-1)/3}\mod{p^2}.$$
\endpro
Proof. Taking $a=-\f 13$ in Theorem 3.1 and noting that
$\ap=\f{p-1}3$, $t=-\f 13$ and
$\b{-1/3}k\b{-2/3}k=\b{2k}k\b{3k}k27^{-k}$ yields the result.

\par\q
 \pro{Theorem 3.4} Let $p$ be a prime with $p\e 1\mod 6$. Then
$$\sum_{k=0}^{p-1}\f{\b{3k}k\b{6k}{3k}}{432^k}\cdot\f p{6k+1}
\e \f{5-2^p}3(-1)^{\f{p-1}6}\b{(p-1)/3}{(p-1)/6} \mod {p^2}.$$
\endpro
Proof.  By (2.2) and (2.4), $H_{[\f p3]}\e -\f 32q_p(3)\mod p$ and
$H_{[\f p6]}\e -2q_p(2)-\f 32q_p(3)\mod p.$  Thus, putting $a=-\f
16$ in Theorem 3.1 and noting that $\ap=\f{p-1}6$ and $t=-\f 16$
gives
$$
\align &\sum_{k=0}^{p-1}\f{\b{3k}k\b{6k}{3k}}{432^k}\cdot\f p{6k+1}
=\f 16\sum_{k=0}^{p-1}\b{-\f 16}k\b{-\f 56}k\f p{k-(-\f 16)}
\\&\e (-1)^{\f{p-1}6}\b{(p-1)/3}{(p-1)/6}\big(1-\f p3(H_{[\f
p3]}-H_{[\f p6]})\big)\\&\e
(-1)^{\f{p-1}6}\b{(p-1)/3}{(p-1)/6}\big(1-\f 23pq_p(2)\big)\mod
{p^2}.
\endalign$$ This proves the theorem.
\section*{4. Congruences for
$\sum_{k=0}^{p-1}\b{2k}k^2\f p{16^k(k+b)}\mod {p^2}$}
\par Let $p$ be an odd prime and $b\in\Bbb Z_p$.
From (2.10) we have a congruence for
$\sum_{k=0}^{p-1}\f{\b{2k}k^2}{16^k(k+b)}$ $\mod {p^2}$ under the
condition $1\le \bp\le \f{p-3}2$. Now we present a congruence for
$\sum_{k=0}^{p-1}\b{2k}k^2\f p{16^k(k+b)}\mod {p^2}$ under the
condition $\bp>\f p2$.
 \pro{Theorem 4.1} Let $p$ be an odd prime,
$b\in\Bbb Z_p$, $\bp>\f p2$ and $s=(b-\bp)/p\not\e -1\mod p$. Then
$$\align
&\sum_{k=0}^{p-1}\f{\b{2k}k^2}{16^k}\cdot\f
p{k+b}\e\sum_{k=0}^{(p-1)/2}\f{\b{2k}k^2}{16^k}\cdot\f p{k+b}\\&\e
\f{\b{2(\bp-\f{p+1}2)}{\bp-\f{p+1}2}^2}{(s+1)16^{\bp-\f{p+1}2}}
\big(1+p(2s+1)(H_{p-\bp} -H_{\bp-\f{p+1}2})\big)
\\&\e \f 1{s+1}\b{\f{p-1}2}{p-\bp}^2
\Big(1+p\big(2q_p(2)+(2s+2)H_{p-\bp}-(2s+1)H_{\bp-\f{p+1}2}\big)\Big)
\\&\e \f{\b{2(p-\bp)}{p-\bp}^2}{(s+1)16^{p-\bp}}
\Big(1+p(2s+2)\big(H_{p-\bp}-H_{\bp-\f{p+1}2}\big)\Big)
 \mod {p^2}.\endalign$$
 \endpro
 Proof.  From [G, (3.100)] we know that
$$\sum_{k=0}^n\b{2k}k(-1)^k\b{n+k}{2k}\f {1}{k+b}
=(-1)^n\f{(b-1)(b-2)\cdots(b-n)}{b(b+1)\cdots(b+n)}.\tag 4.1$$
 From [S3, Lemma 2.2],
$$\b{\f{p-1}2+k}{2k}\e\f{\b{2k}k}{(-16)^k}\Big(1-p^2\sum_{i=1}^k\f
1{(2i-1)^2}\Big) \mod{p^4}\qtq{for}k=1,2,\ldots,\f{p-1}2.$$ Thus,
$$\f{\b{2k}k}{(-16)^k}\e \Big(1+p^2\sum_{i=1}^k\f
1{(2i-1)^2}\Big)\b{\f{p-1}2+k}{2k}
\mod{p^4}\qtq{for}k=1,2,\ldots,\f{p-1}2.\tag 4.2$$ Appealing to
(4.1) and (4.2),
$$\align &\sum_{k=0}^{\f{p-1}2}\f{\b{2k}k^2}{16^k}\cdot\f p{k+b}
\\&\e \sum_{k=0}^{\f{p-1}2}\b{2k}k(-1)^k\b{\f{p-1}2+k}{2k}
\Big(1+p^2\sum_{i=1}^k\f 1{(2i-1)^2}\Big)\f p{k+b}
\\&\e \sum_{k=0}^{\f{p-1}2}\b{2k}k(-1)^k\b{\f{p-1}2+k}{2k}
\f
p{k+b}\\&\qq+p^2\b{2(p-\bp)}{p-\bp}(-1)^{p-\bp}\b{\f{p-1}2+p-\bp}{2(p-\bp)}
\Big(\sum_{i=1}^{p-\bp}\f 1{(2i-1)^2}\Big)\f p{p-\bp+b}
\\&\e  (-1)^{\f{p-1}2}\f{(b-1)(b-2)\cdots(b-\f{p-1}2)}
{b(b+1)\cdots(b+p-\bp-1)\f{b+p-\bp}p(b+p-\bp+1)\cdots (b+\f{p-1}2)}
\\&\qq
+\f
{(-1)^{p-\bp}}{s+1}\b{2(p-\bp)}{p-\bp}\Big(\f{\b{2(p-\bp)}{p-\bp}}
{(-16)^{p-\bp}}-\b{\f{p-1}2+p-\bp}{2(p-\bp)}\Big)
\mod{p^3}.\endalign$$ That is,
$$\aligned \sum_{k=0}^{\f{p-1}2}\f{\b{2k}k^2}{16^k}\cdot\f p{k+b}
&\e \f{(-1)^{p-\bp}\b{b-1}{\f{p-1}2}\b{\f{p-1}2}{\bp-\f{p+1}2}}
{(s+1)\b{-1+p+b-\bp}{p-\bp}\b{-1-(p+b-\bp)}{\bp-\f{p+1}2}}
\\&\q+\f
{(-1)^{p-\bp}}{s+1}\b{2(p-\bp)}{p-\bp}\Big(\f{\b{2(p-\bp)}{p-\bp}}
{(-16)^{p-\bp}}-\b{\f{p-1}2+p-\bp}{2(p-\bp)}\Big) \mod {p^3}.
\endaligned\tag 4.3$$
From (4.3), (2.5) and (4.2),
$$\align&\sum_{k=0}^{(p-1)/2}\f{\b{2k}k^2}{16^k}\cdot\f p{k+b}
\\&\e \f{(-1)^{p-\bp}\b{b-1}{\f{p-1}2}\b{\f{p-1}2}{\bp-\f{p+1}2}}
{(s+1)\b{-1+p+b-\bp}{p-\bp}\b{-1-(p+b-\bp)}{\bp-\f{p+1}2}}
=\f{(-1)^{p-\bp}\b{\bp-1+ps}{\f{p-1}2}\b{\f{p-1}2}{\bp-\f{p+1}2}}
{(s+1)\b{-1+p(s+1)}{p-\bp}\b{-1-p(s+1)}{\bp-\f{p+1}2}}
\\&\e \f{(-1)^{p-\bp}}{s+1}\cdot\f{\b{\bp-1}{\f{p-1}2}
(1+ps\sum_{i=1}^{\f{p-1}2}\f 1{\bp-i})\b{\f{p-1}2}{\bp-\f{p+1}2}}
{\b{-1}{p-\bp}(1+p(s+1)\sum_{i=1}^{p-\bp}\f 1{-i})
\b{-1}{\bp-\f{p+1}2}(1-p(s+1)\sum_{i=1}^{\bp-\f{p+1}2}\f 1{-i})}
\\&=\f 1{s+1}(-1)^{\bp-\f{p+1}2}
\f{\b{\bp-1}{\f{p-1}2}(1+ps(H_{\bp-1}-H_{\bp-\f{p+1}2})\b{\f{p-1}2}{\bp-\f{p+1}2}}{(1-p(s+1)H_{p-\bp})
(1+p(s+1)H_{\bp-\f{p+1}2})}
\\&\e \f 1{s+1}(-1)^{\bp-\f{p+1}2}
\b{\bp-1}{\f{p-1}2}\b{\f{p-1}2}{\bp-\f{p+1}2}
\big(1+ps(H_{p-\bp}-H_{\bp-\f{p+1}2})\big)
\\&\qq\times\big(1+p(s+1)H_{p-\bp}\big) \big(1-p(s+1)H_{\bp-\f{p+1}2}\big)
\\&\e (-1)^{\bp-\f{p+1}2}
\f{\b{\f{p-1}2+\bp-\f{p+1}2}{\bp-\f{p+1}2}\b{\f{p-1}2}{\bp-\f{p+1}2}}{s+1}
\big(1+p(2s+1)(H_{p-\bp} -H_{\bp-\f{p+1}2})\big)
\\&=(-1)^{\bp-\f{p+1}2}\f{\b{2(\bp-\f{p+1}2)}{\bp-\f{p+1}2}\b{\f{p-1}2
+\bp-\f{p+1}2}{2(\bp-\f{p+1}2)}}{s+1} \big(1+p(2s+1)(H_{p-\bp}
-H_{\bp-\f{p+1}2})\big)\\& \e
\f{\b{2(\bp-\f{p+1}2)}{\bp-\f{p+1}2}^2}{(s+1)16^{\bp-\f{p+1}2}}
\big(1+p(2s+1)(H_{p-\bp} -H_{\bp-\f{p+1}2})\big)  \mod
{p^2}.\endalign$$ By [S3, Lemma 2.4], for $k=1,2,\ldots,\f{p-1}2$,
$$\b{\f{p-1}2}k\e \f{\b{2k}k}{(-4)^k}\Big(1-p\sum_{i=1}^k\f
1{2i-1}\Big)=\f{\b{2k}k}{(-4)^k}\Big(1-p\big(H_{2k}-\f
12H_k\big)\Big)\mod {p^2}.\tag 4.4$$ Thus,
$$\b{\f{p-1}2}k^2\e \f{\b{2k}k^2}{16^k}\Big(1-p\big(H_{2k}-\f
12H_k\big)\Big)^2\e \f{\b{2k}k^2}{16^k}(1-2pH_{2k}+pH_k)\mod
{p^2}\tag 4.5$$ and so
$$\f{\b{2k}k^2}{16^k}\e
\b{\f{p-1}2}k^2(1+2pH_{2k}-pH_k)\mod{p^2}.\tag 4.6$$ Now, from the
above and (2.3) we deduce that
$$\align&\sum_{k=0}^{(p-1)/2}\f{\b{2k}k^2}{16^k}\cdot\f p{k+b}
\\& \e \f 1{s+1}\b{\f{p-1}2}{\bp-\f{p+1}2}^2
\big(1+2pH_{2\bp-p-1}-pH_{\bp-\f{p+1}2}\big) \\&\qq\times
\big(1+p(2s+1)(H_{p-\bp} -H_{\bp-\f{p+1}2})\big)
 \\&\e\f
1{s+1}\b{\f{p-1}2}{p-\bp}^2 \Big(1+p\big(2H_{2(p-\bp)}
+(2s+1)H_{p-\bp}-(2s+2)H_{\bp-\f{p+1}2}\big)\Big)
\\&\e \f
1{s+1}\b{\f{p-1}2}{p-\bp}^2
\Big(1+p\big(2q_p(2)+(2s+2)H_{p-\bp}-(2s+1)H_{\bp-\f{p+1}2}\big)\Big)
\\&\e \f 1{(s+1)16^{p-\bp}}\b{2(p-\bp)}{p-\bp}^2
\big(1-2pH_{2p-2\bp}+pH_{p-\bp}\big)
\\&\qq\times\Big(1+p\big(2H_{2(p-\bp)}+(2s+1)
H_{p-\bp}-(2s+2)H_{\bp-\f{p+1}2}\big)\Big)
\\&\e
\f{\b{2(p-\bp)}{p-\bp}^2}{(s+1)16^{p-\bp}}
\Big(1+p(2s+2)\big(H_{p-\bp}-H_{\bp-\f{p+1}2}\big)\Big)
 \mod {p^2}.\endalign$$
To complete the proof, we note that  $p\mid \b{2k}k$ for $\f
p2<k<p$.
\par\q
\pro{Corollary 4.1} Let $p$ be an odd prime, $b\in\Bbb Z_p$, $\bp>\f
p2$ and $s=(b-\bp)/p\not\e -1,-\f 12\mod p$. Then
$$\sum_{k=0}^{p-1}\f{\b{2k}k^2}{16^k}\cdot\f p{k+\f
12-b}\e-\f{2s+2}{2s+1} \sum_{k=0}^{p-1}\f{\b{2k}k^2}{16^k}\cdot\f
p{k+b}\mod {p^2}.$$
\endpro
Proof. Set $b'=\f 12-b$ and $s'=(b'-\langle b'\rangle_p)/p$. Then
$\langle b'\rangle_p=p+\f{p+1}2-\bp>\f p2$ and $s'=-\f 32-s$. By
Theorem 4.1,
$$\align&(s'+1)\sum_{k=0}^{p-1}\f{\b{2k}k^2}{16^k}\cdot\f p{k+b'}
\\&\e \b{\f{p-1}2}{p-\langle b'\rangle_p}^2
\Big(1+p\big(2q_p(2)+(2s'+2)H_{p-\langle
b'\rangle_p}-(2s+1)H_{\langle b'\rangle_p-\f{p+1}2}\big)\Big)
\\&= \b{\f{p-1}2}{p-\langle b\rangle_p}^2
\Big(1+p\big(2q_p(2)+(2s+2)H_{p-\langle b\rangle_p}-(2s+1)H_{\langle
b\rangle_p-\f{p+1}2}\big)\Big)
\\&\e (s+1)\sum_{k=0}^{p-1}\f{\b{2k}k^2}{16^k}\cdot\f p{k+b}
\mod {p^2}.
\endalign$$
This yields the result.
\par\q
 \pro{Theorem 4.2} Let $p$ be a prime with $p\e 1\mod 3$ and so
 $p=x^2+3y^2$. Then
$$\sum_{k=0}^{p-1}\b{2k}k^2\f 1{16^k}\cdot\f p{3k+1}
\e\sum_{k=0}^{p-1}\b{2k}k^2\f 1{16^k}\cdot\f p{6k+1} \e
4x^2-2p\mod{p^2}.$$
\endpro
Proof. Putting $b=\f 13$ in Theorem 4.1 and noting that
$\bp=\f{2p+1}3$ and $s=-\f 23$ gives
$$\align &\sum_{k=0}^{p-1}\b{2k}k^2\f
1{16^k}\cdot\f p{3k+1}\\&\e\b{\f{p-1}2}{\f{p-1}3}^2
\Big(1+p\big(2q_p(2)+\f 23H_{\f{p-1}3}+\f 13H_{\f{p-1}6}\big)\Big)
\\&\e\b{\f{p-1}2}{\f{p-1}3}^2
\Big(1+p\Big(2q_p(2)+\f 23\Big(-\f 32q_p(3)\Big)+\f
13\Big(-2q_p(2)-\f 32q_p(3)\Big)\Big)\Big)
\\&=\b{\f{p-1}2}{\f{p-1}6}^2\Big(1+p\Big(\f 43q_p(2)-\f
32q_p(3)\Big)\Big) \mod {p^2}.\endalign$$
By [Y] or [BEW, Theorem
9.4.4],
$$\b{\f{p-1}2}{\f{p-1}6}\e \Big(2x-\f p{2x}\Big)
\Big(-1+\f 23pq_p(2)-\f 34pq_p(3)\Big)\mod {p^2}.$$ Thus,
$$\b{\f{p-1}2}{\f{p-1}6}^2\e (4x^2-2p)\Big(1-\f 43pq_p(2)+\f
32pq_p(3)\Big)\mod{p^2}.$$ Hence
$$\align &\sum_{k=0}^{p-1}\b{2k}k^2\f
1{16^k}\cdot\f p{3k+1}\\&\e (4x^2-2p)\Big(1-p\Big(\f 43q_p(2)-\f
32q_p(3)\Big)\Big)\Big(1+p\Big(\f 43q_p(2)-\f 32q_p(3)\Big)\Big) \e
4x^2-2p\mod {p^2}.\endalign$$ Taking $b=\f 13$ in Corollary 4.1
gives
 $$\sum_{k=0}^{p-1}\b{2k}k^2\f 1{16^k}\cdot\f p{6k+1} \e
 \sum_{k=0}^{p-1} \b{2k}k^2\f 1{16^k}\cdot\f p{3k+1} \mod {p^2}.$$
Thus, the theorem is proved.
\par\q
\pro{Theorem 4.3} Let $p$ be a prime with $p>3$, $p\e 1,3\mod 8$ and
so $p=x^2+2y^2$. Then
$$\align \big(2-(-1)^{\f{p-1}2}\big)\sum_{k=0}^{p-1}\b{2k}k^2\f 1{16^k}\cdot\f
p{8k+1}&\e \big(2+(-1)^{\f{p-1}2}\big)\sum_{k=0}^{p-1}\b{2k}k^2\f
1{16^k}\cdot\f p{8k+3}\\&\e 4x^2-2p\mod {p^2}.\endalign$$
\endpro
Proof. Since $\f 38=\f 12-\f 18$, from Corollary 4.1 we only need to
prove the congruence for $\sum_{k=0}^{p-1}\b{2k}k^2\f 1{16^k}\cdot\f
p{8k+1}\mod {p^2}.$ Set $b=\f 18$. For $p\e 1\mod 8$,
$\bp=\f{7p+1}8>\f p2$, $p-\bp=\f{p-1}8$, $\bp-\f{p+1}2=\f{3(p-1)}8$
and $s=\f{b-\bp}p =-\f 78$. By Theorem 4.1,
$$\f 18\sum_{k=0}^{p-1}\b{2k}k^2\f 1{16^k}\cdot\f
p{k+1/8}\e \b{\f{p-1}2}{[p/8]}^2 \Big(1+p\Big(2q_p(2)+\f
34H_{[\f{3p}8]}+\f 14H_{[\f p8]}\Big)\Big)\mod {p^2}.$$
 For $p\e 3\mod 8$,
$\bp=\f{5p+1}8>\f p2$, $p-\bp=[\f{3p}8]$, $\bp-\f{p+1}2=[\f p8]$ and
$s=\f{b-\bp}p =-\f 58$. By Theorem 4.1,
$$\f 38\sum_{k=0}^{p-1}\b{2k}k^2\f 1{16^k}\cdot\f
p{k+1/8}\e \b{\f{p-1}2}{[p/8]}^2 \Big(1+p\Big(2q_p(2)+\f
34H_{[\f{3p}8]}+\f 14H_{[\f p8]}\Big)\Big)\mod {p^2}.$$ By (2.3),
$H_{[\f p8]}+H_{[\f{3p}8]}\e 2H_{[\f p4]}-2q_p(2)\e -8q_p(2)\mod p$.
Hence, for $p\e 1,3\mod 8$,
$$\align \big(2-(-1)^{\f{p-1}2}\big)\sum_{k=0}^{p-1}\b{2k}k^2\f 1{16^k}\cdot\f
p{8k+1}&\e \b{\f{p-1}2}{[\f p8]}^2 \Big(1+p\Big(2q_p(2)+\f
34(-8q_p(2))-\f 24H_{[\f p8]}\Big)\Big)
\\&=\b{\f{p-1}2}{[\f p8]}^2 \Big(1-p\Big(4q_p(2)+\f 12H_{[\f p8]}\Big)\Big)
\mod {p^2}.\endalign$$ On the other hand, taking $n=[\f p8]$ and
$t=-\f 14$ in [S7, Theorem 3.2(ii)] gives
$$\align \sum_{k=0}^{p-1}\f{\b{2k}k\b{4k}{2k}}{128^k}
&=\sum_{k=0}^{p-1}\b{-\f 14}k\b{-\f 34}k\f 1{2^k} \e \b{\f{p-1}2}{[
p/8]}\Big(1+p\Big(\f 34H_{[\f p4]}-\f 14H_{[\f p8]}+\f
14q_p(2)\Big)\Big)
\\&\e \b{\f{p-1}2}{[p/8]}\Big(1-p\Big(2q_p(2)+\f 14H_{[\f
p8]}\Big)\Big)\mod {p^2}.\endalign$$ From [S5, Theorem 4.3],
$$\sum_{k=0}^{p-1}\f{\b{2k}k\b{4k}{2k}}{128^k}\e
(-1)^{[\f p8]+\f{p-1}2}\Big(2x-\f p{2x}\Big)\mod {p^2},$$ where $x\e
1\mod 4$. Thus,
$$\align 4x^2-2p&\e \Big(2x-\f p{2x}\Big)^2\e
\Big(\sum_{k=0}^{p-1}\f{\b{2k}k\b{4k}{2k}}{128^k}\Big)^2
\e\b{\f{p-1}2}{[p/8]}^2\Big(1-p\Big(2q_p(2)+\f 14H_{[\f
p8]}\Big)\Big)^2 \\&\e\b{\f{p-1}2}{[p/8]}^2\Big(1-2p\Big(2q_p(2)+\f
14H_{[\f p8]}\Big)\Big)\\&\e
\big(2-(-1)^{\f{p-1}2}\big)\sum_{k=0}^{p-1}\b{2k}k^2\f
1{16^k}\cdot\f p{8k+1} \mod {p^2}.\endalign$$ This proves the
theorem.
\par\q
\pro{Theorem 4.4} Let $p$ be a prime, $p\e 1\mod 4$ and so
$p=x^2+y^2$ with $2\nmid x$. Then
$$\align \Big(3\Ls p3-2\Big)\sum_{k=0}^{p-1}\b{2k}k^2\f 1{16^k}\cdot\f
p{12k+1}&\e\Big(3\Ls p3+2\Big)\sum_{k=0}^{p-1}\b{2k}k^2\f
1{16^k}\cdot\f p{12k+5}
\\&\e 4x^2-2p\mod {p^2}.\endalign$$
\endpro
Proof. Since $\f 5{12}=\f 12-\f 1{12}$, from Corollary 4.1 we only
need to prove the congruence for $\sum_{k=0}^{p-1}\b{2k}k^2\f
1{16^k}\cdot\f p{12k+1}\mod {p^2}.$ Set $b=\f 1{12}$. For $p\e 1\mod
{12}$, $\bp=\f{11p+1}{12}>\f p2$, $p-\bp=\f{p-1}{12}$,
$\bp-\f{p+1}2=\f{5(p-1)}{12}$ and $s=\f{b-\bp}p =-\f {11}{12}$. By
Theorem 4.1,
$$\sum_{k=0}^{p-1}\b{2k}k^2\f 1{16^k}\cdot\f
p{12k+1}\e\b{\f{p-1}2}{[\f p{12}]}^2\Big(1+p\Big(2q_p(2)+\f 16H_{[\f
p{12}]}+\f 56H_{[\f{5p}{12}]}\Big)\Big)\mod {p^2}.$$ For $p\e 5\mod
{12}$, $\bp=\f{7p+1}{12}>\f p2$, $p-\bp=\f{5p-1}{12}$,
$\bp-\f{p+1}2=\f{p-5}{12}$ and $s=\f{b-\bp}p =-\f {7}{12}$. By
Theorem 4.1,
$$\Big(1-\f 7{12}\Big)\sum_{k=0}^{p-1}\b{2k}k^2\f 1{16^k}\cdot\f
p{k+1/12}\e\b{\f{p-1}2}{[\f p{12}]}^2\Big(1+p\Big(2q_p(2)+\f
16H_{[\f p{12}]}+\f 56H_{[\f{5p}{12}]}\Big)\Big)\mod {p^2}.$$ Since
$[\f p{12}]+[\f {5p}{12}]=\f{p-1}2$, from (2.3) and (2.4),
$$H_{[\f{5p}{12}]}+H_{[\f p{12}]}\e 2H_{[\f p6]}-2q_p(2)
\e -6q_p(2)-3q_p(3)\mod p.\tag 4.7$$ Thus,
$$\align 2q_p(2)+\f 16H_{[\f
p{12}]}+\f 56H_{[\f{5p}{12}]} &\e2q_p(2) -\f 46H_{[\f p{12}]}+\f
56(-6q_p(2)-3q_p(3))
\\&=-3q_p(2)-\f 52q_p(3)
-\f 23H_{[\f p{12}]}\mod {p^2}.\endalign$$ On the other hand,
putting $n=[\f p{12}]$ and $t=-\f 16$ in [S7, Theorem 3.1(ii)] we
deduce that
$$\sum_{k=0}^{p-1}\f{\b{3k}k\b{6k}{3k}}{864^k}=\sum_{k=0}^{p-1}
\b{-\f 16}k\b{-\f 56}k\f 1{2^k}\e \b{\f{p-1}2}{[\f
p{12}]}\Big(1+p\Big(\f 56H_{[\f p6]}-\f 13H_{[\f p{12}]}+\f
16q_p(2)\Big)\Big)\mod {p^2}.$$ By [S6, Theorem 3.2],
$\sum_{k=0}^{p-1}\f{\b{3k}k\b{6k}{3k}}{864^k}\e 2c-\f p{2c}\mod
{p^2}$, where $$c=\cases x&\t{if $p\e 1\mod {12}$ and $3\nmid x$,}
\\-x&\t{if $p\e 1\mod {12}$ and $3\mid x$,}
\\y&\t{if $p\e 5\mod {12}$ and $y\e x\mod 3$}
\endcases\tag 4.8$$ and $x\e 1\mod 4$. By (2.4), $H_{[\f p6]}\e -2q_p(2)-\f 32q_p(3)\mod p$. Thus,
$$2c-\f p{2c}\e \b{\f{p-1}2}{[\f p{12}]}\Big(1-p\Big(\f 32q_p(2)
+\f 54q_p(3)+\f 13H_{[\f p{12}]}\Big)\Big)\mod {p^2}\tag 4.9$$ and
so
$$\align 4c^2-2p&\e \b{\f{p-1}2}{[\f p{12}]}^2\Big(1-p\Big(3q_p(2)
+\f 52q_p(3)+\f 23H_{[\f p{12}]}\Big)\Big)\mod {p^2}\\& \e \cases
\sum_{k=0}^{p-1}\b{2k}k^2\f 1{16^k}\cdot\f p{12k+1}\mod {p^2}&\t{if
$p\e 1\mod {12}$,} \\5\sum_{k=0}^{p-1}\b{2k}k^2\f 1{16^k}\cdot\f
p{12k+1}\mod {p^2}&\t{if $p\e 5\mod {12}$.}\endcases\endalign$$ To
complete the proof we note that $4c^2-2p=4y^2-2p\e -(4x^2-2p)\mod
{p^2}$ for $p\e 5\mod{12}$.
\par\q
\par Based on calculations with Maple, we pose the following
conjectures.
 \pro{Conjecture 4.1} Let $p$ be a prime such that  $p\e
1\mod 6$ and so $p=x^2+3y^2$. Then
$$\align \sum_{k=0}^{p-1}\b{2k}k^2\f 1{16^k}\cdot\f
p{3k+1}\e \sum_{k=0}^{p-1}\b{2k}k^2\f 1{16^k}\cdot\f p{6k+1} \e
4x^2-2p-\f{p^2}{4x^2}\mod{p^3}.\endalign$$
\endpro
\pro{Conjecture 4.2} Let $p$ be a prime with $p>3$, $p\e 1,3\mod 8$
and so $p=x^2+2y^2$. Then
$$\align \big(2-(-1)^{\f{p-1}2}\big)\sum_{k=0}^{p-1}\b{2k}k^2\f 1{16^k}\cdot\f
p{8k+1}&\e \big(2+(-1)^{\f{p-1}2}\big)\sum_{k=0}^{p-1}\b{2k}k^2\f
1{16^k}\cdot\f p{8k+3}\\&\e 4x^2-2p-\f{p^2}{4x^2}\mod
{p^3}.\endalign$$
\endpro

\pro{Conjecture 4.3} Let $p$ be a prime, $p\e 1\mod 4$ and so
$p=x^2+y^2$ with $2\nmid x$. Then
$$\align \Big(3\Ls p3-2\Big)\sum_{k=0}^{p-1}\b{2k}k^2\f 1{16^k}\cdot\f
p{12k+1}&\e\Big(3\Ls p3+2\Big)\sum_{k=0}^{p-1}\b{2k}k^2\f
1{16^k}\cdot\f p{12k+5}
\\&\e 4x^2-2p-\f{p^2}{4x^2}\mod {p^3}.\endalign$$
\endpro

\section*{5. Congruences for Domb and Almkvist-Zudilin numbers}
\par This section is devoted to proving some conjectures on
congruences for Domb and Almkvist-Zudilin numbers.
 \pro{Theorem 5.1} Let $p$ be a prime with $p>3$.
Then
$$ \sum_{n=0}^{p-1}\f{D_n}{16^n}\e \sum_{n=0}^{p-1}\f{D_n}{4^n}
 \e
\cases 4x^2-2p\mod{p^2}&\t{if $p\e 1\mod 3$ and so
$p=x^2+3y^2$,}\\0\mod {p^2}&\t{if $p\e 2\mod 3$.}
\endcases$$
\endpro
Proof. From [S8] we know that
$$D_n=(-1)^n\sum_{k=0}^n\b{2k}k^2\b{3k}k\b{n+2k}{3k}(-16)^{n-k}.\tag 5.1$$
From [G, (1.52)],
$$\sum_{n=r}^m\b nr=\b{m+1}{r+1}.\tag 5.2$$
Thus,
$$\align \sum_{n=0}^{p-1}\f{D_n}{16^n}&=\sum_{n=0}^{p-1}
\sum_{k=0}^n\b{2k}k^2\b{3k}k\b{n+2k}{3k}\f 1{(-16)^k}
\\&=\sum_{k=0}^{p-1}\b{2k}k^2\b{3k}k\f
1{(-16)^k}\sum_{n=k}^{p-1}\b{n+2k}{3k}
=\sum_{k=0}^{p-1}\b{2k}k^2\b{3k}k\f
1{(-16)^k}\b{p+2k}{3k+1}.\endalign$$ If $\f{2p}3<k<p$, then
$\b{2k}k^2\b{3k}k=\f{(2k)!(3k)!}{k!^5}\e 0\mod {p^3}$. If $\f
p2<k<\f{2p}3$, then $3k+2\le 2p<p+2k$ and so
$$\b{p+2k}{3k+1}=\b{p+2k}{p-1-k}=
\f{(p+2k)(p+2k-1)\cdots(3k+2)}{(p-1-k)!}\e 0\mod p.$$ Hence $p^3\mid
\b{2k}k^2\b{p+2k}{3k+1}$ and therefore
$$ \aligned\sum_{n=0}^{p-1}\f{D_n}{16^n}&=\sum_{k=0}^{p-1}\b{2k}k^2\b{3k}k\f
1{(-16)^k}\b{p+2k}{3k+1}\\&\e \sum_{k=0}^{(p-1)/2}\b{2k}k^2\b{3k}k\f
1{(-16)^k}\b{p+2k}{3k+1}\mod{p^3}.\endaligned\tag 5.3$$ For $1\le
k<\f p2$ with $k\not=\f{p-1}3$,
$$\align \b{3k}k\b{p+2k}{3k+1}&=
\f{(p+2k)\cdots(p+1)p(p-1)\cdots(p-k)}{k!(2k)!(3k+1)}
\\&=\f
p{3k+1}\cdot\f{(p^2-1^2)\cdots(p^2-k^2)(p+k+1)\cdots(p+2k)}{k!(2k)!}
\\&\e \f
p{3k+1}\cdot\f{(-1^2)\cdots(-k^2)(k+1)\cdots
2k}{k!(2k)!}\Big(1+p\Big(\f 1{k+1}+\cdots+\f 1{2k}\Big)\Big)
\\&=\f{(-1)^k}{3k+1}(p+p^2(H_{2k}-H_k))\mod{p^3}.\endalign$$
For $k=\f{p-1}3$ we see that
$$\align
&\b{3k}k\b{p+2k}{3k+1}=\b{p-1}{\f{2(p-1)}3}\b{p+\f{2(p-1)}3}p
\\&=\f{(p-1)(p-2)\cdots(p-\f{2(p-1)}3)}{\f{2(p-1)}3!}\cdot
\f{(p+\f{2(p-1)}3)\cdots(p+1)p(p-1)\cdots(\f{2(p-1)}3+1)} {p!}
\\&=\f{(p^2-1^2)(p^2-2^2)\cdots(p^2-\sls{2(p-1)}3^2)}{\f{2(p-1)}3!^2}
\e 1-p^2\sum_{i=1}^{2(p-1)/3}\f 1{i^2}\mod {p^4}.
\endalign$$
Let $\{U_n\}$ be given by
$$ U_0=1,\q U_n=-2\sum_{k=1}^{[n/2]}\b n{2k}U_{n-2k}\q(n\ge 1).\tag 5.4$$
Using [S1, Theorem 5.1(a)] and [S4, Theorem 3.3],
$$ \sum_{i=1}^{2(p-1)/3}\f 1{i^2}\e \sum_{i=1}^{2(p-1)/3}
\f 1{(p-i)^2}=\sum_{x=\f{p-1}3+1}^{p-1}\f 1{x^2} =\sum_{x=1}^{p-1}\f
1{x^2}-\sum_{x=1}^{(p-1)/3}\f 1{x^2}  \e 0-3U_{p-3}\mod p.$$ From
[S4, Theorem 3.2(iii)],
$$\sum_{i=\f{p-1}3+1}^{2(p-1)/3}\f 1i=\sum_{i=1}^{2(p-1)/3}
\f 1i-\sum_{i=1}^{(p-1)/3}\f 1i
=\sum_{i=1}^{2(p-1)/3}\f{(-1)^{i-1}}i\e 3pU_{p-3}\mod{p^2}.$$ Hence
for $p\e 1\mod 3$,
$$\b{p-1}{\f{p-1}3}\b{p+\f{2(p-1)}3}p
\e 1-p^2\sum_{i=1}^{2(p-1)/3}\f 1{i^2} \e 1+3p^2U_{p-3}\e
1+p(H_{\f{2(p-1)}3}-H_{\f{p-1}3})\mod {p^3}.$$
 Therefore, from (5.3)
and the above congruence for $\b{3k}k\b{p+2k}{3k+1}\mod {p^3}$ we
derive that
$$\sum_{n=0}^{p-1}\f{D_n}{16^n}
\e\sum_{k=0}^{(p-1)/2}\b{2k}k^2\f 1{16^k}\cdot\f
p{3k+1}(1+p(H_{2k}-H_k))\mod{p^3}.\tag 5.5$$
 Note that $H_{\f{2(p-1)}3}=H_{p-1-\f{p-1}3}\e H_{\f{p-1}3}\mod
p$. We get $\sum_{n=0}^{p-1}\f{D_n}{16^n}
\e\sum_{k=0}^{(p-1)/2}\b{2k}k^2\f 1{16^k}\cdot\f p{3k+1}\mod{p^2}.$
Since
$$\b{p-1}{\f{p-1}3}\b{p+\f{p-1}3}p=\f{(p^2-1^2)(p^2-2^2)\cdots
(p^2-\sls{p-1}3^2)}{\sls{p-1}3!^2}\e 1\mod{p^2},$$ from [S10, p.137]
we know that $\sum_{n=0}^{p-1}\f{D_n}{4^n}
\e\sum_{k=0}^{(p-1)/2}\b{2k}k^2\f 1{16^k}\cdot\f p{3k+1}\mod{p^2}.$
Now applying (2.11) and Theorem 4.2 yields the result.
\par\q
\par {\bf Remark 5.1}  Suppose that $p$ is a prime with $p>3$.
In [S10] the author proved that $\sum_{n=0}^{p-1}\f{D_n}{4^n}\e
0\mod{p^2}$ for $p\e 2\mod 3$.

\par\q Let $\{b_n\}$ be Almkvist-Zudilin numbers. From [S8, (5.1) and (5.2)] we know that
$$\aligned
b_n&=\sum_{k=0}^{[n/3]}\b{2k}k^2\b{4k}{2k}\b{n+k}{4k}(-3)^{n-3k}
\\&=\sum_{k=0}^n\b{2k}k^2\b{4k}{2k}\b{n+3k}{4k}(-27)^{n-k}.\endaligned
\tag 5.6$$

 \pro{Lemma 5.1} Suppose that $p$ is a prime with $p>3$. Then
$$\sum_{n=0}^{p-1}\f{b_n}{(-3)^n}
\e\sum_{n=0}^{p-1}\f{b_n}{(-27)^n}
\e\sum_{k=0}^{p-1}\f{\b{2k}k\b{3k}k}{27^k} \cdot\f p{4k+1}\mod
{p^2}.$$
\endpro
Proof. Using (5.6) and (5.2),
$$\align \sum_{n=0}^{p-1}\f{b_n}{(-3)^n}
&=\sum_{n=0}^{p-1}\sum_{k=0}^{[n/3]}\b{2k}k^2\b{4k}{2k}\b{n+k}{4k}
\f 1{(-27)^k}
\\&=\sum_{k=0}^{[p/3]}\b{2k}k^2\b{4k}{2k}\f
1{(-27)^k}\sum_{n=3k}^{p-1}\b{n+k}{4k}
\\&=\sum_{k=0}^{[p/3]}\b{2k}k^2\b{4k}{2k}\f
1{(-27)^k}\b{p+k}{4k+1}
\\&=\sum_{k=0}^{[p/3]}\f{(p+k)(p+k-1)\cdots(p+1)p(p-1)\cdots(p-3k)}{(4k+1)\cdot
(-27^k)\cdot k!^4}
\\&=\sum_{k=0}^{[p/3]}\f p{4k+1}\cdot
\f{(p^2-1^2)(p^2-2^2)\cdots(p^2-k^2)\cdot
(p-(k+1))\cdots(p-3k)}{k!^4\cdot (-27)^k}.\endalign$$ Thus,
$$\align &\sum_{n=0}^{p-1}\f{b_n}{(-3)^n}
\\&\e \sum_{k=0}^{[p/3]}\f p{(4k+1)(-27)^k\cdot k!^4}\cdot
(-1^2)(-2^2)\cdots(-k^2)\Big(1-p^2\sum_{i=1}^k\f 1{i^2}\Big)
\\&\qq\times(-(k+1))\cdots(-3k)\Big(1-p\sum_{i=k+1}^{3k}\f 1i+p^2\sum_{k+1\le
i<j\le 3k}\f 1{ij}\Big)
\\&=\sum_{k=0}^{[p/3]}\f p{4k+1}\cdot\f{(3k)!}{k!^3\cdot 27^k}
\Big(1-p^2\sum_{i=1}^k\f 1{i^2}\Big) \Big(1-p\sum_{i=k+1}^{3k}\f
1i+\f{p^2}2\Big(\Big(\sum_{i=k+1}^{3k}\f
1i\Big)^2-\sum_{i=k+1}^{3k}\f 1{i^2}\Big)\Big)
\\&\e \sum_{k=0}^{[p/3]}\f p{4k+1}\cdot\f{(3k)!}{k!^3\cdot 27^k}
\Big(1-p\sum_{i=k+1}^{3k}\f 1i+\f{p^2}2\Big(\Big(\sum_{i=k+1}^{3k}\f
1i\Big)^2-\sum_{i=1}^{3k}\f 1{i^2}-\sum_{i=1}^k\f 1{i^2}\Big)\Big)
 \mod{p^3}.\endalign$$
For $k=\f{p-1}4$ we see that $\sum_{i=k+1}^{3k}\f
1i=H_{p-1-\f{p-1}4}-H_{\f{p-1}4}\e 0\mod p$ and
$$\sum_{i=1}^{3(p-1)/4}\f 1{i^2}+\sum_{i=1}^{(p-1)/4}\f 1{i^2}
\e \sum_{k=(p-1)/4+1}^{p-1}\f 1{k^2}+\sum_{k=1}^{(p-1)/4}\f 1{k^2}
=\sum_{k=1}^{p-1}\f 1{k^2}\e 0\mod p.$$ Thus,
$$\sum_{n=0}^{p-1}\f{b_n}{(-3)^n}\e
\sum_{k=0}^{[p/3]}\f p{4k+1}\cdot\f{(3k)!}{k!^3\cdot 27^k}
\Big(1-p\sum_{i=k+1}^{3k}\f 1i\Big)\mod {p^3}\tag 5.7$$ and so
$\sum_{n=0}^{p-1}\f{b_n}{(-3)^n}\e \sum_{k=0}^{[p/3]}\f
p{4k+1}\cdot\f{(3k)!}{k!^3\cdot 27^k}\mod {p^2}.$
 Similarly,
$$\align \sum_{n=0}^{p-1}\f{b_n}{(-27)^n}
&=\sum_{n=0}^{p-1}\sum_{k=0}^{n}\b{2k}k^2\b{4k}{2k}\b{n+3k}{4k} \f
1{(-27)^k}
\\&=\sum_{k=0}^{p-1}\b{2k}k^2\b{4k}{2k}\f
1{(-27)^k}\sum_{n=k}^{p-1}\b{n+3k}{4k}
=\sum_{k=0}^{p-1}\b{2k}k^2\b{4k}{2k}\f 1{(-27)^k}\b{p+3k}{4k+1}
\\&=p+\sum_{k=1}^{p-1}\f{(p+3k)(p+3k-1)\cdots(p+1)p(p-1)\cdots(p-k)}{(4k+1)\cdot
(-27^k)\cdot k!^4}
\\&=p+\sum_{k=1}^{p-1}\f p{4k+1}\cdot
\f{(p^2-1^2)(p^2-2^2)\cdots(p^2-k^2)\cdot
(p+(k+1))\cdots(p+3k)}{k!^4\cdot (-27)^k}.\endalign$$ For $\f
p3<k<p$ we see that $\f p{4k+1}\cdot\f{(3k)!}{k!^3\cdot 27^k}\e
0\mod{p^2}$. Recall that $\sum_{i=k+1}^{3k}\f 1i\e 0\mod p$ for
$k=\f{p-1}4$. From the above we then get
$$\align \sum_{n=0}^{p-1}\f{b_n}{(-27)^n}
&\e p+\sum_{k=1}^{p-1}\f p{(4k+1)(-27)^k\cdot k!^4}
(-1^2)(-2^2)\cdots(-k^2)(k+1)\cdots 3k\Big(1+\sum_{i=k+1}^{3k}\f
pi\Big)
\\&=\sum_{k=0}^{p-1}\f p{4k+1}\cdot\f{(3k)!}{k!^3\cdot 27^k}
 \Big(1+\sum_{i=k+1}^{3k}\f pi\Big)
\e \sum_{k=0}^{[p/3]}\f p{4k+1} \cdot \f{(3k)!}{k!^3\cdot
27^k}\Big(1+\sum_{i=k+1}^{3k}\f pi\Big)
\\&\e\sum_{k=0}^{[p/3]}
\f{\b{2k}k\b{3k}k}{27^k} \cdot\f p{4k+1}\mod {p^2}.
\endalign$$
For $\f{2p}3<k<p$ we have $p^2\mid (3k)!$, for $\f p3<k<\f{2p}3$ we
have $p\nmid (4k+1)$ and $p\mid (3k)!$. Thus,
$$\sum_{k=0}^{p-1}
\f{\b{2k}k\b{3k}k}{27^k} \cdot\f p{4k+1}\e\sum_{k=0}^{[p/3]}
\f{\b{2k}k\b{3k}k}{27^k} \cdot\f p{4k+1}\mod {p^2}.$$ Combining the
above proves the lemma.
\par\q
\pro{Theorem 5.2} Suppose that $p$ is a prime with $p>3$ and $p\e
3\mod 4$. Then
$$\sum_{n=0}^{p-1}\f{b_n}{(-3)^n}\e \sum_{n=0}^{p-1}\f{b_n}{(-27)^n}
\e 0\mod{p^2}.$$
\endpro
Proof. From Lemma 5.1 and (2.12) (with $b=\f 14$),
$$\sum_{n=0}^{p-1}\f{b_n}{(-3)^n}\e \sum_{n=0}^{p-1}\f{b_n}{(-27)^n}
\e  \sum_{k=0}^{p-1}\f{\b{2k}k\b{3k}k}{27^k}\cdot \f p{4k+1} \e
0\mod {p^2}.$$
\par\q
  \pro{Theorem 5.3} Let $p$ be a prime with $p\e 1\mod 4$. Then
$$\sum_{n=0}^{p-1}\f{b_n}{(-3)^n}\e \sum_{n=0}^{p-1}\f{b_n}{(-27)^n}
\e\cases 4x^2-2p\mod {p^2}&\t{if $p\e 1\mod{12}$ and so
$p=x^2+9y^2$,}
\\2p-2x^2\mod {p^2}&\t{if $p\e 5\mod{12}$ and so
$2p=x^2+9y^2$.}
\endcases$$
\endpro
Proof. Set $S(a)=\sum_{k=0}^{p-1}\b ak\b{-1-a}k\f p{4k+1}.$ By
(1.2), $S(-\f 13)=\sum_{k=0}^{p-1}\f{(3k)!}{k!^3\cdot 27^k}\cdot \f
p{4k+1}.$  Applying Lemma 5.1,
$$\sum_{n=0}^{p-1}\f{b_n}{(-3)^n}\e \sum_{n=0}^{p-1}\f{b_n}{(-27)^n}
\e S\Big(-\f 13\Big)=S\Big(-\f 23\Big)\mod {p^2}.\tag 5.8$$ For $p\e
1\mod{12}$, taking $a=-\f 13$, $b=\f 14$ and $m=\f {p-1}{12}-1$ in
Lemma 2.1 yields
$$S\Big(-\f 13\Big)\e \f{\b{-7/12}{(p-1)/12}}{\b{-1/12}{(p-1)/12}}
S\Big(-\f 14-\f p{12}\Big) \mod {p^3}.\tag 5.9$$ For $p\e
5\mod{12}$, taking $a=-\f 23$, $b=\f 14$ and $m=\f {p-5}{12}-1$ in
Lemma 2.1 yields
$$S\Big(-\f 23\Big)\e \f{\b{-11/12}{(p-5)/12}}{\b{-5/12}{(p-5)/12}}
S\Big(-\f 14-\f p{12}\Big) \mod {p^3}.\tag 5.10$$
 Note that
$\sum_{k=1}^{(p-1)/2}\f 1{2k-1}=H_{p-1}-\f 12H_{\f{p-1}2}\e
q_p(2)\mod p$. By [S2, Theorem 3.1(ii)],
$$\sum_{i=0}^{(p-1)/4-1}\f 1{4i+1}=\sum\Sb k=1\\k\e p\mod
4\endSb^{p-1}\f 1k\e \f 34q_p(2)\mod p.\tag 5.11$$ Thus,
$$\sum_{i=0}^{\f{p-1}4-1}\f 1{4i+3}=\sum_{k=1}^{\f{p-1}2}\f 1{2k-1}
-\sum_{i=0}^{\f{p-1}4-1}\f 1{4i+1}\e q_p(2)-\f 34q_p(2)=\f
14q_p(2)\mod p.\tag 5.12$$
  Using (2.5), (5.11) and (5.12),
$$\align S\Big(-\f 14-\f
p{12}\Big)&=\sum_{k=0}^{p-1}\b{-\f 14-\f p{12}}k\b{-\f 34+\f
p{12}}k\f p{4k+1}
\\&\e \sum_{k=0}^{p-1}\b{-\f 14}k\b{-\f 34}k\f p{4k+1}
\Big(1+\sum_{i=0}^{k-1}\f{-\f p{12}}{-\f 14-i}\Big)
\Big(1+\sum_{i=0}^{k-1}\f{\f p{12}}{-\f 34-i}\Big)
\\&= \sum_{k=0}^{p-1}\b{-\f 14}k\b{-\f 34}k\f p{4k+1}
\Big(1+\f 13\sum_{i=0}^{k-1}\f p{4i+1}\Big)\Big(1-\f
13\sum_{i=0}^{k-1}\f p{4i+3}\Big)
\\&\e \sum_{k=0}^{p-1}\b{-\f 14}k\b{-\f 34}k\f p{4k+1}
\Big(1+\f 13\sum_{i=0}^{k-1}\Big(\f p{4i+1}-\f p{4i+3}\Big)\Big)
\\&\e \sum_{k=0}^{p-1}\b{-\f 14}k\b{-\f 34}k\f p{4k+1}
+\f p3\b{-\f 14}{\f{p-1}4}\b{-\f
34}{\f{p-1}4}\sum_{i=0}^{(p-1)/4-1}\Big(\f 1{4i+1}-\f 1{4i+3}\Big)
\\&\e S\Big(-\f 14\Big)+\f p3\cdot
\f{\b{\f{p-1}2}{\f{p-1}4}\b{p-1}{\f{p-1}2}}{64^{\f{p-1}4}}\Big(\f
34q_p(2)-\f 14q_p(2)\Big)\\&\e S\Big(-\f 14\Big)+\f
p6\b{\f{p-1}2}{\f{p-1}4} (-1)^{\f{p-1}4}q_p(2) \mod {p^2}.
\endalign$$
 Now applying Theorem 3.2 gives
$$ S\Big(-\f 14-\f{p}{12}\Big)
\e S\Big(-\f 14\Big)+\f p6\b{\f{p-1}2}{\f{p-1}4}
(-1)^{\f{p-1}4}q_p(2) \e (-1)^{\f{p-1}4}\b{\f{p-1}2}{\f{p-1}4}
\Big(1-\f 13pq_p(2)\Big) \mod {p^2}.$$ Suppose $p=A^2+B^2$ with $A\e
1\mod 4$. From [CDE] or [BEW],
$$\b{(p-1)/2}{(p-1)/4}\e \Big(1+\f 12pq_p(2)\Big)\Big(2A-\f
p{2A}\Big)\mod {p^2}.\tag 5.13$$ Thus,
$$\align S\Big(-\f 14-\f{p}{12}\Big)&\e (-1)^{\f{p-1}4}
\Big(1+\f 12pq_p(2)\Big)\Big(2A-\f p{2A}\Big)\Big(1-\f
13pq_p(2)\Big)
\\&\e (-1)^{\f{p-1}4}\Big(1+\f 16pq_p(2)\Big)\Big(2A-\f
p{2A}\Big) \mod {p^2}.\endalign$$ Combining this with (5.8)-(5.10)
yields
$$\aligned &\sum_{n=0}^{p-1}\f{b_n}{(-3)^n}\e
\sum_{n=0}^{p-1}\f{b_n}{(-27)^n} \\&\e \cases
\f{\b{-7/12}{(p-1)/12}}{\b{-1/12}{(p-1)/12}}
(-1)^{\f{p-1}4}\big(1+\f 16pq_p(2)\big)\big(2A-\f p{2A}\big) \mod
{p^2}&\t{if $p\e 1\mod{12}$,}
\\\f{\b{-11/12}{(p-5)/12}}{\b{-5/12}{(p-5)/12}}
(-1)^{\f{p-1}4}\big(1+\f 16pq_p(2)\big)\big(2A-\f p{2A}\big) \mod
{p^2}&\t{if $p\e 5\mod{12}$.}
\endcases\endaligned\tag 5.14$$
For $p\e 1\mod {12}$ suppose $p=x^2+9y^2$. Using (2.5) we see that
$$\align\f{\b{-\f 7{12}}{\f{p-1}{12}}}{\b{-\f 1{12}}{\f{p-1}{12}}}
&=\f{\b{\f{7(p-1)}{12}-\f{7p}{12}}{\f{p-1}{12}}}
{\b{\f{p-1}{12}-\f{p}{12}}{\f{p-1}{12}}} \e\f{\b{\f{7(p-1)}{12}}
{\f{p-1}{12}}\big(1-\f 7{12}p\sum_{i=0}^{\f{p-1}{12}-1}\f
1{\f{7(p-1)}{12}-i}\big)} {\b{\f{p-1}{12}} {\f{p-1}{12}}\big(1-\f
p{12}\sum_{i=0}^{\f{p-1}{12}-1}\f 1{\f{p-1}{12}-i}\big)}
\\&=\b{\f{7(p-1)}{12}}
{\f{p-1}{12}}\f{1-\f 7{12}p(H_{\f{7(p-1)}{12}}-H_{\f{p-1}2})}{1-\f
p{12}H_{\f{p-1}{12}}}
\\&\e \b{\f{7(p-1)}{12}}
{\f{p-1}{12}}\Big(1-\f 7{12}p(H_{\f{7(p-1)}{12}}-H_{\f{p-1}2})\Big)
\Big(1+\f p{12}H_{\f{p-1}{12}}\Big)
\\&\e \b{\f{7(p-1)}{12}}
{\f{p-1}{12}}\Big(1-\f 7{12}p\big(H_{\f{7(p-1)}{12}}+2q_p(2)\big)+\f
p{12}H_{\f{p-1}{12}}\Big)\mod {p^2}.
\endalign$$
By (4.7), $H_{\f{7(p-1)}{12}}\e H_{\f{5(p-1)}{12}}\e
-6q_p(2)-3q_p(3)-H_{\f{p-1}{12}}\mod p$. By [S3, Lemma 2.5],
$$\align \b{\f{7(p-1)}{12}}
{\f{p-1}{12}}&\e (-1)^{\f{p-1}{12}}\b{\f{p-1}2}{\f{p-1}{12}}
\Big(1+2p\Big(H_{\f{p-1}6}-\f 12H_{\f{p-1}{12}}\Big)\Big) \\& \e
(-1)^{\f{p-1}{12}}\b{\f{p-1}2}{\f{p-1}{12}}
\Big(1-p\big(4q_p(2)+3q_p(3)+H_{\f{p-1}{12}}\big)\Big)\mod {p^2}.
\endalign$$
Therefore,
$$\align
&(-1)^{\f{p-1}{4}}\Big(1+\f 16pq_p(2)\Big)\b{-\f
7{12}}{\f{p-1}{12}}\b{-\f 1{12}}{\f{p-1}{12}}^{-1}\\ &\e \Big(1+\f
16pq_p(2)\Big)\b{\f{p-1}2}{\f{p-1}{12}}
\Big(1-p\big(4q_p(2)+3q_p(3)+H_{\f{p-1}{12}}\big)\Big)
\\&\qq\times\Big(1-\f
7{12}p\big(-4q_p(2)-3q_p(3)-H_{\f{p-1}{12}}\big)+\f
p{12}H_{\f{p-1}{12}}\Big)
\\&\e\b{\f{p-1}2}{\f{p-1}{12}}
\Big(1-p\Big(\f 32q_p(2)+\f 54q_p(3)+\f 13H_{\f{p-1}{12}}\Big)\Big)
\mod {p^2}.
\endalign$$
This together with (5.14) and (4.9) yields
$$\align \sum_{n=0}^{p-1}\f{b_n}{(-3)^n}&\e
\sum_{n=0}^{p-1}\f{b_n}{(-27)^n} \e \b{\f{p-1}2}{\f{p-1}{12}}
\Big(1-p\Big(\f 32q_p(2)+\f 54q_p(3)+\f 13H_{\f{p-1}{12}}\Big)\Big)
\Big(2A-\f p{2A}\Big)
\\&\e \Big(2c-\f p{2c}\Big)\Big(2A-\f p{2A}\Big)
\e\cases 4A^2-2p=4x^2-2p \mod {p^2}&\t{if $3\nmid A$,}
\\-(4A^2-2p)\e 4x^2-2p\mod {p^2}&\t{if $3\mid A$},
\endcases\endalign$$
where $c=A$ or $-A$ according as $3\nmid A$ or $3\mid A$.
\par Now assume $p\e 5\mod {12}$ and $p=A^2+B^2$ with $A\e 1\mod 4$
 and $B\e A\mod 3$. Appealing to (2.5),
$$\align\f{\b{-\f {11}{12}}{\f{p-5}{12}}}{\b{-\f 5{12}}{\f{p-5}{12}}}
&=\f{\b{\f{7p-11}{12}-\f{7p}{12}}{\f{p-5}{12}}}
{\b{\f{p-5}{12}-\f{p}{12}}{\f{p-5}{12}}} \e\f{\b{\f{7p-11}{12}}
{\f{p-5}{12}}\big(1-\f {7}{12}p\sum_{i=0}^{\f{p-5}{12}-1}\f
1{\f{7p-11}{12}-i}\big)} {\b{\f{p-5}{12}} {\f{p-5}{12}}\big(1-\f
p{12}\sum_{i=0}^{\f{p-5}{12}-1}\f 1{\f{p-5}{12}-i}\big)}
\\&=\b{\f{7p-11}{12}}
{\f{p-5}{12}}\f{1-\f {7}{12}p(H_{\f{7p-11}{12}}-H_{\f{p-1}2})}{1-\f
p{12}H_{\f{p-5}{12}}}
\\&\e \b{\f{7p-11}{12}}
{\f{p-5}{12}}\Big(1-\f {7}{12}p(H_{\f{7p-11}{12}}-H_{\f{p-1}2})\Big)
\Big(1+\f p{12}H_{\f{p-5}{12}}\Big)
\\&\e \b{\f{7p-11}{12}}
{\f{p-5}{12}}\Big(1-\f
{7}{12}p\big(H_{\f{7p-11}{12}}+2q_p(2)\big)+\f
p{12}H_{\f{p-5}{12}}\Big)\mod {p^2}.
\endalign$$
By (4.7), $$H_{\f{7p-11}{12}}=H_{p-1-\f{5p-1}{12}}\e
H_{\f{5p-1}{12}} \e -6q_p(2)-3q_p(3)-H_{\f{p-5}{12}}.$$ By [S3,
Lemma 2.5],
$$\align \b{\f{7p-11}{12}}
{\f{p-5}{12}}&\e (-1)^{\f{p-5}{12}}
\b{\f{p-1}2}{\f{p-5}{12}}\Big(1+2p\Big(H_{\f{p-5}6}-\f
12H_{\f{p-5}{12}}\Big)\Big) \\&\e (-1)^{\f{p-5}{12}}
\b{\f{p-1}2}{\f{p-5}{12}}\Big(1-p\Big(4q_p(2)+3q_p(3)+
H_{\f{p-5}{12}}\Big)\Big)\mod {p^2}.\endalign$$ Hence
$$\align\f{\b{-\f {11}{12}}{\f{p-5}{12}}}{\b{-\f 5{12}}{\f{p-5}{12}}}
&\e(-1)^{\f{p-5}{12}}
\b{\f{p-1}2}{\f{p-5}{12}}\Big(1-p\Big(4q_p(2)+3q_p(3)+
H_{\f{p-5}{12}}\Big)\Big)
\\&\qq\times\Big(1-\f
{7}{12}p\big(-6q_p(2)-3q_p(3)-H_{\f{p-5}{12}}+2q_p(2)\big)+\f
p{12}H_{\f{p-5}{12}}\Big)
\\&\e(-1)^{\f{p-5}{12}}
\b{\f{p-1}2}{\f{p-5}{12}}\Big(1-p\Big(\f 53q_p(2)+\f 54q_p(3)+\f 13
H_{\f{p-5}{12}}\Big)\Big)\mod {p^2}
\endalign$$
and so
$$\align &\b{\f{7p-11}{12}}
{\f{p-5}{12}} (-1)^{\f{p-1}4}\Big(1+\f 16pq_p(2)\Big)
\\&\e-\b{\f{p-1}2}{\f{p-5}{12}}\Big(1+\f 16pq_p(2)\Big)
\Big(1-p\Big(\f 53q_p(2)+\f 54q_p(3)+\f 13
H_{\f{p-5}{12}}\Big)\Big)
\\&\e-\b{\f{p-1}2}{\f{p-5}{12}}
\Big(1-p\Big(\f 32q_p(2)+\f 54q_p(3)+\f 13 H_{\f{p-5}{12}}\Big)\Big)
\e -\Big(2B-\f p{2B}\Big)\mod {p^2}
\endalign$$ by (4.9). This together with (5.14) yields
$$\sum_{n=0}^{p-1}\f{b_n}{(-3)^n}\e
\sum_{n=0}^{p-1}\f{b_n}{(-27)^n} \e -\Big(2A-\f p{2A}\Big)
\Big(2B-\f p{2B}\Big)\e -4AB\mod {p^2}.$$ Suppose $2p=x^2+9y^2$.
Then clearly $\sls{x+3y}2^2+\sls{x-3y}2^2=p$ and $\f{x+3y}2\e
\f{x-3y}2\mod 3$. Thus,
$$2p-2x^2=9y^2-x^2=-4\cdot\f{x+3y}2\cdot\f{x-3y}2=-4AB.$$
Thus the theorem is proved.
\par\q

\pro{Corollary 5.1} Let $p$ be a prime of the form $12k+5$. Then
$p=A^2+B^2$ with $A\e 1\mod 4$ and $B\e A\mod 3$, and $2p=x^2+9y^2$.
We have
$$\sum_{k=0}^{p-1}
\f{\b{2k}k^2\b{4k}{2k}}{(-12288)^k}\e -4AB\e -2x^2\mod p.$$
\endpro
Proof. Taking $u=-1/3$ in [S8, Theorem 5.1] gives $\sum_{k=0}^{p-1}
\f{\b{2k}k^2\b{4k}{2k}}{(-12288)^k}\e\sum_{n=0}^{p-1}\f{b_n}{(-3)^n}$
$\mod p.$ Thus the result follows from Theorem 5.3 and its proof.
\par\q

\pro{Theorem 5.4} Let $p>3$ be a prime. Then
$$\align &\sum_{n=0}^{p-1}\f{W_n}{(-3)^n}\e\sum_{k=0}^{[p/3]}\f{\b{2k}k\b{3k}k}{27^k}
\cdot\f p{3k+1}\\&\e \cases -L+\f pL\mod{p^2}&\t{if $3\mid p-1$ and
so $4p=L^2+27M^2$ with $3\mid L-1$,}
\\ -\f p3\big(\f{p-2}3!\big)^3
\mod{p^2}&\t{if $p\e 2\mod 3$.}
\endcases\endalign$$
\endpro
Proof. Using (1.8) and (5.2),
$$\align\sum_{n=0}^{p-1}\f{W_n}{(-3)^n}&
=\sum_{n=0}^{p-1}\sum_{k=0}^{[n/3]}\b{2k}k\b{3k}k\b n{3k}\f
1{(-27)^k}
\\&=\sum_{k=0}^{[p/3]}\b{2k}k\b{3k}k\f 1{(-27)^k}\sum_{n=3k}^{p-1}
\b n{3k}=\sum_{k=0}^{[p/3]}\b{2k}k\b{3k}k\f 1{(-27)^k}\b p{3k+1}
\\&=\sum_{k=0}^{[p/3]}\b{2k}k\b{3k}k\f 1{(-27)^k}\cdot\f
p{3k+1}\cdot\f{(p-1)(p-2)\cdots(p-3k)}{(3k)!}\\& \e
\sum_{k=0}^{[p/3]}\f{\b{2k}k\b{3k}k}{27^k} \cdot\f p{3k+1}
\mod{p^2}.\endalign$$  For $\f p3<k<p$ with $k\not=\f{2p-1}3$ we
have $p\nmid 3k+1$ and $\b{2k}k\b{3k}k=\f{(3k)!}{k!^3}\e 0\mod p$.
Thus,
$$\align &\sum_{k=0}^{p-1}\f{\b{2k}k\b{3k}k}{27^k}\cdot\f p{3k+1}
- \sum_{k=0}^{[p/3]}\f{\b{2k}k\b{3k}k}{27^k}\cdot\f p{3k+1}
\\&\e\cases 0\mod{p^2}&\t{if $p\e 1\mod 3$,}
\\\f{(2p-1)!}{(\f{2p-1}3)!^3}\cdot
\f 1{2\cdot 27^{(2p-1)/3}}\mod{p^2}&\t{if $p\e 2\mod 3$.}
\endcases\endalign$$
For $p\e 1\mod 3$ it is well known that $4p=L^2+27M^2$ with
$L,M\in\Bbb Z$ and $L\e 1\mod 3$. By [BEW, Theorem 9.4.4],
$\b{\f{2(p-1)}3}{\f{p-1}3}\e -L+\f pL\mod{p^2}$. Hence applying the
above and Theorem 3.3,
$$\sum_{n=0}^{p-1}\f{W_n}{(-3)^n}\e
\sum_{k=0}^{p-1}\f{\b{2k}k\b{3k}k}{27^k}\cdot\f p{3k+1} \e
\b{\f{2(p-1)}3}{\f{p-1}3}\e -L+\f pL\mod{p^2}.$$
 For $p\e 2\mod 3$,
$$\align\f{(2p-1)!}{(\f{2p-1}3)!^3}\cdot \f 1{2\cdot
27^{(2p-1)/3}}&
=\f{p(p^2-1^2)(p^2-2^2)\cdots(p^2-(p-1)^2)\cdot\sls{p-2}3!^3}
{\sls{2p-1}3!^3\sls{p-2}3!^3\cdot 2\cdot 3^{2p-1}}
\\&\e p\cdot\Big(\f{p-2}3!\Big)^3\f{\b{p-1}{(p-2)/3}^3}{(p-1)!\cdot 2\cdot 3}
\e \f p6\Big(\f{p-2}3!\Big)^3\mod{p^2}
\endalign$$
and
$$\b{\f{2(p-2)}3}{\f{p-2}3}^{-1}=\f{\big(\f{p-2}3!\big)^2(p-1)(p-2)\cdots(p-\f{p+1}3)}
{(p-1)!}\e -\f 13\big(\f{p-2}3!\big)^3\mod p.$$ From (2.16) and the
above,
$$\align \sum_{n=0}^{p-1}\f{W_n}{(-3)^n}&\e
\sum_{k=0}^{p-1}\f{\b{2k}k\b{3k}k}{27^k}\cdot\f p{3k+1}- \f
p6\Big(\f{p-2}3!\Big)^3\e \f p2\b{\f{2(p-2)}3}{\f{p-2}3}^{-1} -\f
p6\Big(\f{p-2}3!\Big)^3\\&\e - \f p6\Big(\f{p-2}3!\Big)^3- \f
p6\Big(\f{p-2}3!\Big)^3=- \f p3\Big(\f{p-2}3!\Big)^3 \mod{p^2}.
\endalign$$
 Thus, the theorem is proved.
\par\q
\par To conclude the paper, we pose two challenging conjectures.
\pro{Conjecture 5.1} Let $p>3$ be a prime. Then
$$\align &\sum_{k=0}^{p-1}\f{\b{2k}k\b{3k}k}{27^k}\cdot\f p{4k+1}
\\&\e \cases 4x^2-2p-\f{p^2}{4x^2}\mod {p^3}&\t{if $12\mid p-1$ and so
$p=x^2+9y^2$,}
\\2p-2x^2+\f{p^2}{2x^2}\mod {p^3}&\t{if $12\mid p-5$ and so
$2p=x^2+9y^2$,} \\ -\f 5{27}p^2\b{[p/3]}{[p/12]}^{-2}\mod
{p^3}&\t{if $12\mid p-7$,}
\\\f 5{54}p^2\b{[p/3]}{[p/12]}^{-2}\mod {p^3}
&\t{if $12\mid p-11$.}\endcases\endalign$$
\endpro
\par\q
\pro{Conjecture 5.2} Suppose that $p$ is a prime of the form $4k+1$.
Then
$$(-1)^{[\f p{12}]}\b{[\f p3]}{[\f p{12}]}\e\cases
2x\mod p&\t{if $p\e 1\mod{12}$ and so $p=x^2+9y^2$ with $3\mid
x-1$,}
\\x\mod p&\t{if $p\e 5\mod{12}$ and so $2p=x^2+9y^2$ with $3\mid x-1$.}
\endcases$$
\endpro

\end{document}